\newtheorem{theorem}{Theorem}

\newtheorem{corollary}[theorem]{Corollary}
\newtheorem{definition}[theorem]{Definition}
\newtheorem{remark}[theorem]{Remark}

\newtheorem{lemma}[theorem]{Lemma}
\newtheorem{proposition}[theorem]{Proposition}

\documentclass[notitlepage,leqno]{article}%
\usepackage{amssymb}
\usepackage{amsfonts}
\usepackage[leqno]{amsmath}
\usepackage{graphicx}
\usepackage[utf8]{inputenc}
\begin{document}

\title{Optimal Dividend Strategies for Two Collaborating Insurance Companies}
\author{Hansj\"{o}rg Albrecher\thanks{Department of Actuarial Science, Faculty of
Business and Economics, University of Lausanne, CH-1015 Lausanne and Swiss
Finance Institute. Supported by the Swiss National Science Foundation Project
$200020\_143889$.}, Pablo Azcue\thanks{Departamento de Matematicas,
Universidad Torcuato Di Tella. Av. Figueroa Alcorta 7350 (C1428BIJ) Ciudad de
Buenos Aires, Argentina.} and Nora Muler$^{\dag}$}
\date{}
\maketitle
\abstract{We consider a two-dimensional optimal dividend problem in the context of two insurance companies
with compound Poisson surplus processes, who collaborate by paying each other's deficit when possible.
We solve the stochastic control problem of maximizing the weighted sum of expected discounted dividend payments
(among all admissible dividend strategies) until ruin of both companies, by extending results of univariate optimal
control theory.  In the case that the dividends paid by the two companies are equally weighted, the value function
of this problem compares favorably with the one of merging the two companies completely. We identify this optimal value function  as the smallest viscosity supersolution of the respective
Hamilton-Jacobi-Bellman equation and provide an iterative approach to approximate it numerically.
Curve strategies are identified as the natural analogue of barrier strategies in this two-dimensional context.
A numerical example is given for which such a curve strategy is indeed optimal among all admissible
dividend strategies, and for which this collaboration mechanism also outperforms the suitably weighted optimal dividend
strategies of the two stand-alone companies.}

\section{Introduction}

Ever since de Finetti \cite{DeFin57} proposed in 1957 to measure the value of
an insurance portfolio by the expected discounted sum of dividends paid during
the lifetime of the portfolio, it has been of particular interest to determine
the optimal dividend payment strategy which maximizes this quantity. More than
that, this field of research over the years turned out to be a challenging and
fascinating area, combining tools from analysis, probability and stochastic
control. In 1969, Gerber \cite{Ger68} showed that if the free surplus of an
insurance portfolio is modelled by a compound Poisson risk model, it is
optimal to pay dividends according to a so-called band strategy, which
collapses to a barrier strategy for exponentially distributed claim amounts.
Whereas Gerber found this result by taking a limit of an associated discrete
problem, this optimal dividend problem was studied with techniques of modern
stochastic control theory in Azcue and Muler \cite{azmu05}, see e.g. Schmidli
\cite{schmidli08} for a detailed overview. Since then the optimal dividend
problem was studied for many different model setups, objective functions and
side constraints (we refer to Albrecher and Thonhauser \cite{AlTho09} and
Avanzi \cite{Avanzi} for surveys on the subject). A barrier strategy with
barrier $b$ pays out dividends whenever the surplus level of the portfolio is
above $b$, so that the surplus level stays at $b$, and pays no dividends below
that barrier $b$. The most general criteria currently available for barrier
strategies to be optimal can be found in Loeffen and Renaud \cite{loeren}. The
optimality of barrier strategies when including the time value of ruin was
studied in \cite{ThAl}, and when including capital injections by shareholders
in Kulenko and Schmidli \cite{KuSch}.\newline

All these control problems have been formulated and studied in the
one-dimensional framework. However, in recent years there has been an
increased interest in risk theory in considering the dynamics of several
connected insurance portfolios simultaneously, see e.g. Asmussen and Albrecher
\cite[Ch.XIII.9]{asal} for an overview. Ruin probability expressions for a
two-dimensional risk process are studied in Avram et al.
\cite{Avram08b,Avram08a} for simultaneous claim arrivals and proportional
claim sizes and recently in Badila et al. \cite{badila} and Ivanovs and Boxma
\cite{IvBox} in a more general framework. In Azcue and Muler \cite{azmu2}, the
problem of optimally transferring capital between two portfolios in the
presence of transaction costs was considered, see also Badescu et al.
\cite{BGL}. Czarna and Palmowski \cite{czpal} study the dividend problem and
impulse control for two insurance companies who share 
claim payments and premiums in some specified proportion for a particular dividend
strategy. It turns out that these multi-dimensional problems, albeit
practically highly relevant, quickly become very intricate and explicit
solutions can typically not be obtained without very strong assumptions.
\newline

In this paper, we would like to extend the optimal dividend problem from
univariate risk theory to a two-dimensional setup of two collaborating
companies. The collaboration consists of paying the deficit ('bailing out') of
the partner company if its surplus is negative and if this financial help can
be afforded with the current own surplus level. We solve the problem of
maximizing the weighted sum of expected discounted dividend payments until
ruin of both companies. A natural question in this context is whether such a
collaboration procedure can be advantageous over merging the two companies
completely; we will show that this is the case when the dividends paid by the
two companies are equally weighted. For criteria of a merger being an
advantage over keeping two stand-alone companies under pre-defined barrier
strategies and marginal diffusion processes, see e.g. Gerber and Shiu
\cite{gershiu06}, for the performance of another pre-defined risk and profit
sharing arrangement, see e.g. Albrecher and Lautscham \cite{AlLau}. Our goal
here is, however, to address the general problem of identifying the optimal
dividend strategy (among all admissible dividend strategies) for each company
under this collaboration framework. This leads to a fully two-dimensional
stochastic control problem, and to the question what the natural analogues of
the optimal univariate barrier strategies are in two dimensions. The
particular structure of the collaboration implemented in this paper will turn
out not to be essential, so the techniques may be applicable to other
risk-sharing mechanisms as well. Yet, the concrete specification allows to
carry through the necessary analysis of the stochastic control problem
explicitly by way of example. \newline

The rest of the paper is organized as follows. In Section \ref{sec2} we
introduce the model and the stochastic control problem in detail and derive
some simple properties of the corresponding value function $V$. In Section
\ref{sec3} we prove that $V$ is a viscosity solution of the corresponding
Hamilton-Jacobi-Bellman equation for independent surplus processes, and in
Section \ref{sec4} we show that $V$ is in fact its smallest viscosity
supersolution. Section \ref{sec5} provides an iterative approach to
approximate the value function $V$, together with the analogous verification
steps at each iteration step. Section \ref{sec6} discusses the stationary
dividend strategies that appear in our model, and in Section \ref{sec7} we establish \emph{curve
strategies} as the appropriate analogues of the univariate barrier strategies.
Finally, Section \ref{sec8} shows how to constructively search for optimal
curve strategies and in Section \ref{sec9} an explicit numerical example for
the symmetric (and equally weighted) case with exponentially distributed claim
sizes is worked out for which such a curve strategy is indeed optimal among
all admissible bivariate dividend strategies. It is then also illustrated that
for this case the proposed type of collaboration is preferable to adding the
best-possible stand-alone profits.

\section{Model}

\label{sec2}

We consider two insurance companies, Company One and Company Two, which have
an agreement to collaborate. Let us call $X_{t}$ the free surplus of Company
One and $Y_{t}$ the one of Company Two. We assume that the free surplus of
each of the companies follows a Cram\'{e}r-Lundberg process, i.e. a compound
Poisson process with drift given by%

\begin{equation}
\left\{
\begin{array}
[c]{l}%
X_{t}=x+p_{1}t-%
%TCIMACRO{\tsum \nolimits_{i=1}^{N_{t}^{1}}}%
%BeginExpansion
{\textstyle\sum\nolimits_{i=1}^{N_{t}^{1}}}
%EndExpansion
U_{i}^{(1)}\\
Y_{t}=y+p_{2}t-%
%TCIMACRO{\tsum \nolimits_{i=1}^{N_{t}^{2}}}%
%BeginExpansion
{\textstyle\sum\nolimits_{i=1}^{N_{t}^{2}}}
%EndExpansion
U_{i}^{(2)},
\end{array}
\right.  \label{freeSurplus}%
\end{equation}
where $x$\ and $y$ are the respective initial surplus levels; $p_{1}$\ and
$p_{2}$\ are the respective premium rates; $U_{i}^{{\small (k)}}\ $is the size
of the $i$-th claim\ of Company $k$, which are i.i.d. random variables with
continuous distribution $F^{k}$ for $k=1,2$; $N_{t}^{1}\ $and $N_{t}^{2}$ are
Poisson processes with intensity $\lambda_{1}\ $and $\lambda_{2}$,
respectively. We assume here that the processes $N_{t}^{1}$, $N_{t}^{2}$\ and
the random variables $U_{i}^{(1)}$, $U_{i}^{(2)}$ are all independent of each
other, and $p_{j}>\lambda_{j} E(U_{i}^{(j)})$, $j=1,2$.

There is a rule of collaboration signed by the two companies: if the current
surplus of Company One becomes negative, Company Two should cover the exact
deficit of Company One as long as it does not ruin itself, and vice versa.
Ruin of a company hence occurs when its surplus becomes negative and the other
company cannot cover this deficit.

A simulated surplus trajectory under this collaboration rule is shown in
Figure 1.1.%

\[%
%TCIMACRO{\FRAME{itbpFU}{4.0465in}{2.5175in}{0in}{\Qcb{Fig. 1.1: Surplus
%process under the collaboration rules.}}{}{trayectoriapaper.eps}%
%width 4.0465in;  height 2.5175in;  depth 0in;  original-width 3.9972in;
%original-height 2.4768in;  cropleft "0";  croptop "1";  cropright "1";
%cropbottom "0";  filename 'TrayectoriaPaper.eps';file-properties "XNPEU";}}}%
%BeginExpansion
{\parbox[b]{4.0465in}{\begin{center}
\includegraphics[
height=2.5175in,
width=4.0465in
]%
{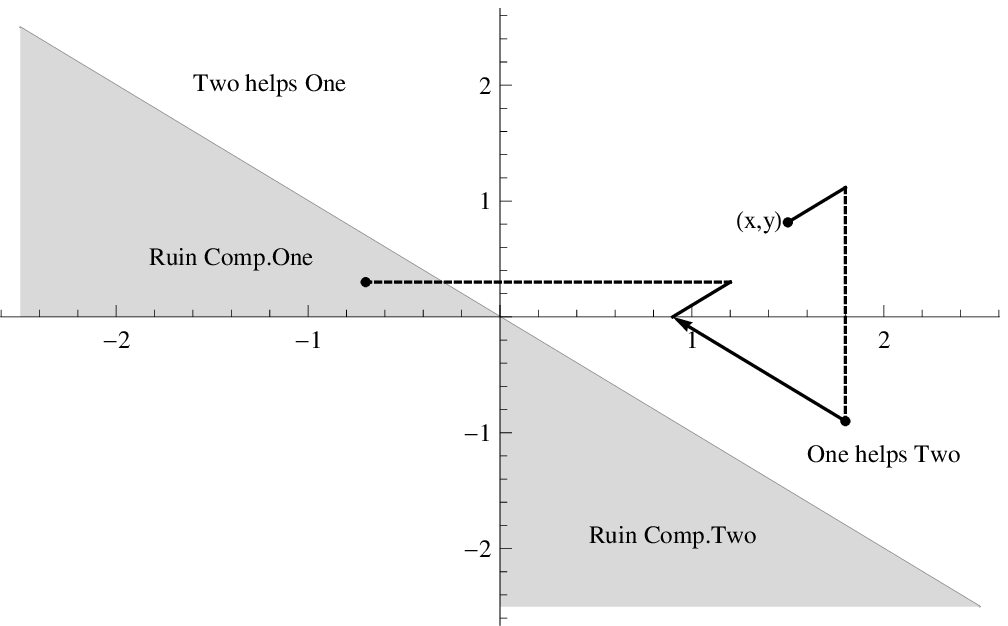}%
\\
Fig. 1.1: Surplus process under the collaboration rules.
\end{center}}}%
%EndExpansion
\]

Both companies use part of their surplus to pay dividends to their
shareholders. The dividend payment strategy $\overline{L}=\left(  L_{t}%
^{1}{\small ,}L_{t}^{2}\right)  $\ is the total amount of dividends paid by
the two companies up to time $t.$ Let us call $\tau_{i}^{k}$ the arrival time
of the $i$-th claim of company $k,$ with $k=1,2$. We define the associated
controlled process $\left(  X_{t}^{\overline{L}},Y_{t}^{\overline{L}}\right)
$ with initial surplus levels $(x,y)$ as%

\begin{equation}
\left\{
\begin{array}
[c]{c}%
X_{t}^{\overline{L}}=X_{t}-L_{t}^{1}+C_{t}^{2,1}-C_{t}^{1,2}\\
Y_{t}^{\overline{L}}=Y_{t}-L_{t}^{2}+C_{t}^{1,2}-C_{t}^{2,1},
\end{array}
\right.  \label{Transferencia}%
\end{equation}
where%

\[
C_{t}^{2,1}=%
%TCIMACRO{\dsum \limits_{i=1}^{N_{t}^{1}}}%
%BeginExpansion
{\displaystyle\sum\limits_{i=1}^{N_{t}^{1}}}
%EndExpansion
I_{\left\{  X_{\tau_{i}^{1}}^{\overline{L}}<0,Y_{\tau_{i}^{1}}^{\overline{L}%
}+X_{\tau_{i}^{1}}^{\overline{L}}\geq0\right\}  }\left\vert X_{\tau_{i}^{1}%
}^{\overline{L}}\right\vert
\]
corresponds to the cumulative amount transferred from Company Two to Company
One up to time $t$ in order to cover the deficit of Company One and%

\[
C_{t}^{1,2}=%
%TCIMACRO{\dsum \limits_{i=1}^{N_{t}^{2}}}%
%BeginExpansion
{\displaystyle\sum\limits_{i=1}^{N_{t}^{2}}}
%EndExpansion
I_{\left\{  Y_{\tau_{i}^{2}}^{\overline{L}}<0,Y_{\tau_{i}^{2}}^{\overline{L}%
}+X_{\tau_{i}^{2}}^{\overline{L}}\geq0\right\}  }\left\vert Y_{\tau_{i}^{2}%
}^{\overline{L}}\right\vert
\]
corresponds to the cumulative amount transferred from Company One to Company
Two up to time $t$ in order to cover the deficit of Company Two.

Let us call $\overline{\tau}$ the time at which only one company remains
(because it cannot cover the deficit of the other), more precisely,%

\begin{equation}
\overline{\tau}=\inf\left\{  t\geq0:X_{t}^{\overline{L}}+Y_{t}^{\overline{L}%
}<0\right\}  \text{.} \label{DefinicionTauRuina}%
\end{equation}

The process $\left(  X_{t}^{\overline{L}},Y_{t}^{\overline{L}}\right)  $ is
defined for $t\leq\overline{\tau}$. We say that the dividend payment strategy
$\overline{L}=\left(  L_{t}^{1}{\small ,}L_{t}^{2}\right)  _{t\leq
\overline{\tau}}$ is \textit{admissible} if it is non-decreasing,
c\`{a}gl\`{a}d (left continuous with right limits), predictable with respect
to the filtration generated by the bivariate process $\left(  X_{t}%
,Y_{t}\right)  $, and satisfies
\[
\left\{
\begin{array}
[c]{c}%
L_{t}^{1}\leq X_{t}+C_{t}^{2,1}-C_{t}^{1,2},\\
L_{t}^{2}\leq Y_{t}+C_{t}^{1,2}-C_{t}^{2,1}\text{.}%
\end{array}
\right.
\]
This last condition means that the companies are not allowed to pay more
dividends than their current surplus. Let us call $\mathbf{R}_{+}^{2}$ the
first quadrant. We denote by $\Pi_{x,y}$ the set of admissible dividend
strategies with initial surplus levels $(x,y)\in\mathbf{R}_{+}^{2}$. Our
objective is to maximize the weighted average of the expected discounted
dividends paid by the two companies until ruin of both companies. Note that
after time $\overline{\tau}$, the surviving company can continue to pay 
dividends up to its own ruin. Let us define $V_{k}^{0}$ ($k=1,2$) as the
optimal value function of the one-dimensional problem of maximizing the
expected discounted dividends until ruin of Company $k$ alone. So, for any
initial surplus levels $(x,y)\in\mathbf{R}_{+}^{2}$, we can write down the
optimal value function as%
\begin{equation}
V(x,y)=\sup_{\overline{L}\in\Pi_{x,y}}V_{\overline{L}}(x,y),
\label{DefinitionV}%
\end{equation}
where%

\begin{equation}
V_{\overline{L}}(x,y)=E_{x,y}\left(  a_{1}\left(
%TCIMACRO{\tint \limits_{0}^{\overline{\tau}}}%
%BeginExpansion
{\textstyle\int\limits_{0}^{\overline{\tau}}}
%EndExpansion
e^{-\delta s}dL_{s}^{1}+e^{-\delta\overline{\tau}}V_{1}^{0}(X_{\overline{\tau
}}^{\overline{L}})\right)  +a_{2}\left(
%TCIMACRO{\tint \limits_{0}^{\overline{\tau}}}%
%BeginExpansion
{\textstyle\int\limits_{0}^{\overline{\tau}}}
%EndExpansion
e^{-\delta s}dL_{s}^{2}+e^{-\delta\overline{\tau}}V_{2}^{0}(Y_{\overline{\tau
}}^{\overline{L}})\right)  \right)  \text{.} \label{Definicion VL}%
\end{equation}
Here, $\delta>0$ is a constant discount factor, and $a_{1}\in\lbrack0,1]$ and
$a_{2}=1-a_{1}$ are the weights of the dividends paid by Company One and
Company Two respectively. The functions $V_{k}^{0}$ ($k=1,2$) are zero in
$\left(  -\infty,0\right)  $ so depending on which company goes to ruin at
$\overline{\tau}$, either $V_{1}^{0}(X_{\overline{\tau}})=0$ or $V_{2}%
^{0}(Y_{\overline{\tau}})=0$. The optimal dividend strategy corresponding to
(\ref{DefinitionV}) may be regarded as the best dividend payment strategy
from the point of view of a shareholder who owns a proportion $ma_{1}$ of the
total shares of Company One and a proportion $ma_{2}$ of the total shares of Company
Two for some $0<m\leq\min\left\{  1/a_{1},1/a_{2}\right\}  $. An important
particular case is $a_{1}=a_{2}=1/2,$ in which the dividends paid by the two
companies are equally weighted (for an earlier example of weighting separate
terms in the objective function in the univariate dividend context, see Radner
and Shepp \cite{radsh}).

\begin{remark}
\label{Comparacion Merger} In case the two companies are owned
by the same shareholders, another possibility of collaboration between the two
companies is merging, in which case the companies put together all their surplus, pay the claims of both companies and pay dividends up to time
$\tau$ at which the joined surplus becomes negative (see e.g. Gerber and Shiu
\cite{gershiu06})). Given the initial surplus levels $(x,y)$, we can interpret
any admissible dividend payment strategy $\left(  L_{t}\right)  _{t\geq0}$ for
the merger as an admissible collaborating one as follows%
\[
L_{t}^{1}=%
%TCIMACRO{\tint \limits_{0}^{t}}%
%BeginExpansion
{\textstyle\int\limits_{0}^{t}}
%EndExpansion
\tfrac{X_{t}^{\overline{L}}}{X_{t}^{\overline{L}}+Y_{t}^{\overline{L}}}%
\,dL_{t}\text{, }L_{t}^{2}=%
%TCIMACRO{\tint \limits_{0}^{t}}%
%BeginExpansion
{\textstyle\int\limits_{0}^{t}}
%EndExpansion
\tfrac{Y_{t}^{\overline{L}}}{X_{t}^{\overline{L}}+Y_{t}^{\overline{L}}}%
\,dL_{t}.
\]

Since $L_{t}$ is constant for $t\geq\tau=\overline{\tau},$ the surviving
company does not pay any dividends here. So $V_{\overline{L}}$ defined in
(\ref{Definicion VL}) for $a_{1}=a_{2}=1/2,$ satisfies
\[
2V_{\overline{L}}(x,y)>2E_{x,y}\left(  \frac{1}{2}%
%TCIMACRO{\tint \limits_{0}^{\overline{\tau}}}%
%BeginExpansion
{\textstyle\int\limits_{0}^{\overline{\tau}}}
%EndExpansion
e^{-\delta s}dL_{s}^{1}+\frac{1}{2}%
%TCIMACRO{\tint \limits_{0}^{\overline{\tau}}}%
%BeginExpansion
{\textstyle\int\limits_{0}^{\overline{\tau}}}
%EndExpansion
e^{-\delta s}dL_{s}^{2}\right)  =E_{x+y}(%
%TCIMACRO{\tint \limits_{0}^{\tau}}%
%BeginExpansion
{\textstyle\int\limits_{0}^{\tau}}
%EndExpansion
e^{-\delta s}dL_{s}).
\]
The last expected value is the value function of the merger dividend strategy
$\left(  L_{t}\right)  _{t\geq0}.$ We conclude that the optimal collaborating
strategy for equally weighted dividend payments is better than the optimal
merger strategy.
\end{remark}

Both optimal value functions $V_{1}^{0}$ and $V_{2}^{0}$ corresponding to
the stand-alone companies have an ultimately linear growth with slope one and
they are Lipschitz, see for instance Azcue and Muler \cite{azmu05}. Let us
state some basic results about regularity and growth at infinity of the
optimal value function $V$ defined in (\ref{DefinitionV}). From now on, let us
call $\lambda:=\lambda_{1}+\lambda_{2}$ and $p:=a_{1}p_{1}+a_{2}p_{2}.$

\begin{lemma}
\label{V_GrowthCondition}The optimal value function is well defined and
satisfies%
\[
a_{1}x+a_{2}y+\frac{p}{\delta+\lambda}\leq V(x,y)\leq a_{1}x+a_{2}y+\frac
{p}{\delta}%
\]
for all $(x,y)\in\mathbf{R}_{+}^{2}.$
\end{lemma}

\textit{Proof.}

Let us first prove the second inequality. Note that $V(x,y)$ increases when
the Poisson intensities $\lambda_{1}$ and $\lambda_{2}$ decrease, but the
optimal value function for the problem with parameters $\lambda_{1}%
=\lambda_{2}=0$ is
\[
a_{1}x+a_{2}y+\int\nolimits_{0}^{\infty}e^{-\delta s}pds=a_{1}x+a_{2}%
y+\frac{p}{\delta},
\]
which corresponds to the value function of the strategy in which each company
pays immediately the initial surplus and then it pays the incoming premium
forever as dividends.

In order to obtain the first inequality, consider the admissible strategy
$\overline{L}_{0}=\left(  L_{t}^{1}{\small ,}L_{t}^{2}\right)  $ in which each
company pays the initial surplus immediately, and then pays the incoming
premium up to time $\overline{\tau}$ that coincides with the first claim
arrival time $\tau_{1}=\tau_{1}^{1}\wedge\tau_{1}^{2}$; we have,%

\[%
\begin{array}
[c]{lll}%
V(x,y) & \geq & V_{\overline{L}_{0}}(x,y)\\
& \geq & a_{1}x+a_{2}y+E_{x,y}\left(
%TCIMACRO{\tint \limits_{0}^{\tau_{1}}}%
%BeginExpansion
{\textstyle\int\limits_{0}^{\tau_{1}}}
%EndExpansion
e^{-\delta s}pds\right) \\
& = & a_{1}x+a_{2}y+\frac{p}{\delta+\lambda}\text{.}
\end{array}
\]
\hfill${\small \square}$

\begin{lemma}
\label{V increasing_LocallyLip}The optimal value function $V$ is increasing,
locally Lipschitz and satisfies for any $(x,y)\in\mathbf{R}_{+}^{2},$
\[
a_{1}h\leq V(x+h,y)-V(x,y)\leq(e^{(\delta+\lambda)h/p_{1}}-1)V(x,y)
\]
and%
\[
a_{2}h\leq V(x,y+h)-V(x,y)\leq(e^{(\delta+\lambda)h/p_{2}}-1)V(x,y)
\]
for any $h>0$.
\end{lemma}

\textit{Proof.}

Let us prove the inequalities at the top, the ones at the bottom are similar.
Given any $\varepsilon>0$, take an admissible strategy $\overline{L}\in
\Pi_{x,y}$ such that $V_{\overline{L}}(x,y)\geq V(x,y)-\varepsilon$. We define
the strategy $\overline{L}^{1}\in\Pi_{x+h,y}$ for $h>0$ as follows: Pay
immediately an amount $h$ of the surplus of Company One as dividends and then
follow the strategy $\overline{L}$. We have that
\[
V_{\overline{L}^{1}}\left(  x+h,y\right)  =V_{\overline{L}}(x,y)+a_{1}h
\]
and so%

\begin{equation}
V\left(  x+h,y\right)  \geq V_{\overline{L}}(x,y)+a_{1}h>V(x,y)+a_{1}%
h-\varepsilon. \label{Desigualdad1}%
\end{equation}
Consider also an admissible strategy $\overline{L}^{2}\in\Pi_{x+h,y}$ such
that $V(x+h,y)\geq V_{\overline{L}^{2}}(x+h,y)-\varepsilon$ and define the
admissible strategy $\overline{L}^{3}\in\Pi_{x,y}$ which, starting with
surplus $\left(  x,y\right)  $ pays no dividends until%

\[
\tilde{\tau}=\inf\left\{  t\geq0:X_{t}^{\overline{L}^{3}}\geq x+h,\text{
}Y_{t}^{\overline{L}^{3}}\geq y\right\}  ,
\]
at time $\tilde{\tau}$ pays either $X_{t}^{\overline{L}^{3}}-(x+h)$ from the
surplus of Company One or $Y_{t}^{\overline{L}^{3}}-y$ from the surplus of
Company Two, depending on which of these differences is positive, and then
follows strategy $\overline{L}^{2}\in\Pi_{x+h,y}$. In the event of no claims,
$\tilde{\tau}=t_{0}:=h/p_{1};$ since the probability of no claims until
$t_{0}$ is $e^{-\lambda t_{0}}$, we get%

\begin{equation}
V(x,y)\geq V_{\overline{L}^{3}}(x,y)\geq V_{\overline{L}^{2}}%
(x+h,y)e^{-(\delta+\lambda)t_{0}}\geq(V(x+h,y)-\varepsilon)e^{-(\delta
+\lambda)t_{0}}. \label{Desigualdad2}%
\end{equation}
From (\ref{Desigualdad1}) and (\ref{Desigualdad2}), we get the inequalities at
the top.\hfill$~{\small \square}$

\section{Hamilton-Jacobi-Bellman equation}

\label{sec3}

In order to obtain the Hamilton-Jacobi-Bellman (HJB) equation associated to
the optimization problem (\ref{DefinitionV}), we need to state the so called
Dynamic Programming Principle (DPP). The proof that this holds is similar to
the one given in Lemma 1.2 of Azcue and Muler \cite{azmu} and uses that $V$ is
increasing and continuous in $\mathbf{R}_{+}^{2}$.

\begin{lemma}
\label{DPP}For any initial surplus $\left(  x,y\right)  $ in $\mathbf{R}%
_{+}^{2}$ and any stopping time $\tau$, we can write%

\[%
\begin{array}
[c]{l}%
V(x,y)\\%
\begin{array}
[c]{ll}%
= & \sup\limits_{\overline{L}\in\Pi_{x,y}}(E_{x,y}(a_{1}%
%TCIMACRO{\tint \limits_{0}^{\tau\wedge\overline{\tau}}}%
%BeginExpansion
{\textstyle\int\limits_{0}^{\tau\wedge\overline{\tau}}}
%EndExpansion
e^{-\delta s}dL_{s}^{1}+a_{2}%
%TCIMACRO{\tint \limits_{0}^{\tau\wedge\overline{\tau}}}%
%BeginExpansion
{\textstyle\int\limits_{0}^{\tau\wedge\overline{\tau}}}
%EndExpansion
e^{-\delta s}dL_{s}^{2}+e^{-\delta(\tau\wedge\overline{\tau})}I_{\{\tau
\wedge\overline{\tau}<\overline{\tau}\}}V(X_{\tau\wedge\overline{\tau}%
}^{\overline{L}},Y_{\tau\wedge\overline{\tau}}^{\overline{L}})\\
& +e^{-\delta(\tau\wedge\overline{\tau})}I_{\{\tau\wedge\overline{\tau
}=\overline{\tau}\}}(a_{1}V_{1}^{0}(X_{\overline{\tau}}^{\overline{L}}%
)+a_{2}V_{2}^{0}(Y_{\overline{\tau}}^{\overline{L}})))).
\end{array}
\end{array}
\]

\end{lemma}

The HJB equation of this optimization problem is%

\begin{equation}
\max\left\{  \mathcal{L}(V)(x,y),a_{1}-V_{x}(x,y),a_{2}-V_{y}(x,y)\right\}
=0, \label{HJB}%
\end{equation}
where%

\begin{equation}
\mathcal{L}(V)(x,y)=V_{x}(x,y)p_{1}+V_{y}(x,y)p_{2}-\left(  \delta
+\lambda\right)  V(x,y)+\mathcal{I}(V)(x,y)+U(x,y), \label{DefinicionL}%
\end{equation}

\begin{equation}%
\begin{array}
[c]{ccc}%
\mathcal{I}(V)(x,y) & = & \lambda_{1}\int_{0}^{x}V(x-\alpha,y)dF^{1}%
(\alpha)+\lambda_{1}\int_{x}^{x+y}V(0,x+y-\alpha)dF^{1}(\alpha)\\
&  & +\lambda_{2}\int_{0}^{y}V(x,y-\alpha)dF^{2}(\alpha)+\lambda_{2}\int
_{y}^{x+y}V(x+y-\alpha,0)dF^{2}(\alpha),
\end{array}
\nonumber
\end{equation}

and%
\begin{equation}
U(x,y)=\lambda_{1}a_{2}V_{2}^{0}(y)(1-F^{1}(x+y))+\lambda_{2}a_{1}V_{1}%
^{0}(x)(1-F^{2}(x+y)). \label{DefinicionU}%
\end{equation}

Since the optimal value function $V$ is locally Lipschitz but possibly not
differentiable at certain points, we cannot say that $V$ is a solution of the
HJB equation, we prove instead that $V$ is a viscosity solution of the
corresponding HJB equation. Let us define this notion (see Crandell and Lions
\cite{cralio} and Soner \cite{soner} for further details).

\begin{definition}
\label{Viscosity}A locally Lipschitz function $\overline{u}:\mathbf{R}_{+}%
^{2}\rightarrow\mathbf{R}$\ is a viscosity supersolution of (\ref{HJB})\ at
$(x,y)\in\mathbf{R}_{+}^{2}$\ if any continuously differentiable function
$\varphi:\mathbf{R}_{+}^{2}\rightarrow\mathbf{R}\ $with $\varphi
(x,y)=\overline{u}(x,y)$ such that $\overline{u}-\varphi$\ reaches the minimum
at $\left(  x,y\right)  $\ satisfies
\[
\max\left\{  \mathcal{L}(\varphi)(x,y),a_{1}-\varphi_{x}(x,y),a_{2}%
-\varphi_{y}(x,y)\right\}  \leq0.\
\]
A function $\underline{u}:$ $\mathbf{R}_{+}^{2}\rightarrow\mathbf{R}$\ is a
viscosity subsolution\ of (\ref{HJB}) at $(x,y)\in\mathbf{R}_{+}^{2}$\ if any
continuously differentiable function $\psi:\mathbf{R}_{+}^{2}\rightarrow
\mathbf{R}\ $with $\psi(x,y)=\underline{u}(x,y)$ such that $\underline{u}%
-\psi$\ reaches the maximum at $\left(  x,y\right)  $ satisfies
\[
\max\left\{  \mathcal{L}(\psi)(x,y),a_{1}-\psi_{x}(x,y),a_{2}-\psi
_{y}(x,y)\right\}  \geq0\text{.}%
\]
A function $u:\mathbf{R}_{+}^{2}\rightarrow\mathbf{R}$ which is both a
supersolution and subsolution at $(x,y)\in\mathbf{R}_{+}^{2}$ is called a
viscosity solution of (\ref{HJB})\ at $\left(  x,y\right)  \in\mathbf{R}%
_{+}^{2}.$
\end{definition}

\begin{proposition}
\label{Prop V is a viscosity supersolution}$V$ is a viscosity supersolution of
the HJB equation (\ref{HJB}) at any $\left(  x,y\right)  $ with $x>0$ and
$y>0$.
\end{proposition}

\textit{Proof. }Given initial surplus levels $x>0$, $y>0$ and any $l_{1}$
$\geq0$, $l_{2}$ $\geq0$, let us consider the admissible strategy
$\overline{L}$ where Company One and Two pay dividends with constant rates
$l_{1}$ and $l_{2}$ respectively and $\overline{\tau}$ is defined as in
(\ref{DefinicionTauRuina}). Let $\varphi$ be a test function for the
supersolution of (\ref{HJB}) at $(x,y)$ with $x>0$ and $y>0.$ As before,
denote $\tau_{1}^{1}$ and $\tau_{1}^{2}$ as the arrival time of the first
claim of Company One and Two respectively, and $\tau_{1}=\tau_{1}^{1}%
\wedge\tau_{1}^{2}$. We have for $t<\tau_{1}$,%

\[
\left\{
\begin{array}
[c]{l}%
X_{t}^{\overline{L}}=x+\left(  p_{1}-l_{1}\right)  t\text{,}\\
Y_{t}^{\overline{L}}=y+\left(  p_{2}-l_{2}\right)  t\text{.}%
\end{array}
\right.
\]
Note that $N_{t}^{1}+N_{t}^{2}$ is a Poisson process with intensity $\lambda$,
because the arrival times of the two companies are independent. We have from
Lemma \ref{DPP} that%

\[%
\begin{array}
[c]{lll}%
\varphi(x,y) & = & V(x,y)\\
& \geq & E_{x,y}(a_{1}\int\nolimits_{0}^{\tau_{1}\wedge t}e^{-\delta
\,s}\,l_{1}ds+a_{2}\int\nolimits_{0}^{\tau_{1}\wedge t}e^{-\delta\,s}%
\,l_{2}ds)\\
&  & +E_{x,y}\left(  e^{-\delta\,(\tau_{1}\wedge t)}I_{\{\tau_{1}\wedge
t<\overline{\tau}\}}V(X_{\tau_{1}\wedge t}^{\overline{L}},Y_{\tau_{1}\wedge
t}^{\overline{L}}))\right) \\
&  & +E_{x,y}\left(  e^{-\delta\overline{\tau}}I_{\{\tau_{1}\wedge
t=\overline{\tau}\}}(a_{1}V_{1}^{0}(X_{\overline{\tau}}^{\overline{L}}%
)+a_{2}V_{2}^{0}(Y_{\overline{\tau}}^{\overline{L}}))\right) \\
& \geq & E_{x,y}(a_{1}\int\nolimits_{0}^{\tau_{1}\wedge t}e^{-\delta
\,s}\,l\,_{1}ds+a_{2}\int\nolimits_{0}^{\tau_{1}\wedge t}e^{-\delta
\,s}\,l\,_{2}ds)\\
&  & +E_{x,y}\left(  e^{-\delta\,(\tau\wedge t)}\varphi(X_{\tau\wedge
t}^{\overline{L}},Y_{\tau\wedge t}^{\overline{L}})I_{\{\tau_{1}\wedge
t<\overline{\tau}\}}\right) \\
&  & +E_{x,y}\left(  e^{-\delta\overline{\tau}}I_{\{\tau_{1}\wedge
t=\overline{\tau}\}}(a_{1}V_{1}^{0}(X_{\overline{\tau}}^{\overline{L}}%
)+a_{2}V_{2}^{0}(Y_{\overline{\tau}}^{\overline{L}}))\right)  .
\end{array}
\]
We can write%
\[%
\begin{array}
[c]{l}%
E_{x,y}\left(  e^{-\delta\,(\tau\wedge t)}\varphi(X_{\tau\wedge t}%
^{\overline{L}},Y_{\tau\wedge t}^{\overline{L}})I_{\{\tau_{1}\wedge
t<\overline{\tau}\}}\right) \\%
\begin{array}
[c]{ll}%
= & E_{x,y}(I_{\{t<\tau_{1}\text{ }\}}e^{-\delta t}\varphi(X_{\tau_{1}\wedge
t}^{\overline{L}},Y_{\tau_{1}\wedge t}^{\overline{L}}))\\
& +E_{x,y}(I_{\{\tau_{1}=\tau_{1}\wedge t<\overline{\tau}\text{and }\tau
_{1}=\tau_{1}^{1}\text{ }\}}e^{-\delta\tau_{1}^{1}}\varphi(X_{\tau_{1}^{1}%
}^{\overline{L}},Y_{\tau_{1}^{1}}^{\overline{L}}))\\
& +E_{x,y}(I_{\{\tau_{1}=\tau_{1}\wedge t<\overline{\tau}\text{and }\tau
_{1}=\tau_{1}^{2}\text{ }\}}e^{-\delta\tau_{1}^{2}}\varphi(X_{\tau_{1}^{2}%
}^{\overline{L}},Y_{\tau_{1}^{2}}^{\overline{L}})).
\end{array}
\end{array}
\]
So we obtain%

\begin{equation}%
\begin{array}
[c]{lll}%
0 & \geq & \lim_{t\rightarrow0^{+}}E_{x,y}\left(  \frac{a_{1}\int
\nolimits_{0}^{\tau_{1}\wedge t}e^{-\delta\,s}\,l\,_{1}ds+a_{2}\int
\nolimits_{0}^{\tau_{1}\wedge t}e^{-\delta\,s}\,l\,_{2}ds}{t}\right) \\
&  & +\lim_{t\rightarrow0^{+}}\frac{e^{-(\lambda+\delta)t}\varphi(x+\left(
p_{1}-l_{1}\right)  t,y+\left(  p_{2}-l_{2}\right)  t)-\varphi(x,y)}{t}\\
&  & +\lim_{t\rightarrow0^{+}}E_{x,y}\left(  \frac{I_{\{\tau_{1}=\tau
_{1}\wedge t<\overline{\tau}\text{and }\tau_{1}=\tau_{1}^{1}\text{ }%
\}}e^{-\delta\tau_{1}^{1}}\varphi(X_{\tau_{1}^{1}}^{\overline{L}},Y_{\tau
_{1}^{1}}^{\overline{L}})}{t}\right) \\
&  & +\lim_{t\rightarrow0^{+}}E_{x,y}\left(  \frac{I_{\{\tau_{1}=\tau
_{1}\wedge t<\overline{\tau}\text{and }\tau_{1}=\tau_{1}^{2}\text{ }%
\}}e^{-\delta\tau_{1}^{2}}\varphi(X_{\tau_{1}^{2}}^{\overline{L}},Y_{\tau
_{1}^{2}}^{\overline{L}})}{t}\right) \\
&  & +\lim_{t\rightarrow0^{+}}E_{x,y}\left(  \frac{e^{-\delta\overline{\tau}%
}I_{\{\tau_{1}\wedge t=\overline{\tau}\}}(a_{1}V_{1}^{0}(X_{\overline{\tau}%
}^{\overline{L}})+a_{2}V_{2}^{0}(Y_{\overline{\tau}}^{\overline{L}}))}%
{t}\right) \\
& = & a_{1}l_{1}+a_{2}l_{2}-\left(  \delta+\lambda\right)  \varphi(x,y)\\
&  & +\left(  p_{1}-l_{1}\right)  \varphi_{x}(x,y)+\left(  p_{2}-l_{2}\right)
\varphi_{y}(x,y)+\mathcal{I}(\varphi)(x,y)\\
&  & +U(x,y)\text{.}%
\end{array}
\nonumber
\end{equation}
Therefore,%

\[
0\geq\mathcal{L}(\varphi)(x,y)+l_{1}\left(  a_{1}-\varphi_{x}(x,y)\right)
+l_{2}\left(  a_{2}-\varphi_{y}(x,y)\right)  .
\]
Taking $l_{1}=l_{2}=0$, $l_{1}\rightarrow\infty$ with $l_{2}=0$, and
$l_{2}\rightarrow\infty$ with $l_{1}=0$, we obtain%

\[
\max\left\{  \mathcal{L}(\varphi)(x,y),a_{1}-\varphi_{x}(x,y),a_{2}%
-\varphi_{y}(x,y)\right\}  \leq0.\hfill{\small \square}%
\]

\begin{proposition}
\label{Prop V Subsolution}$V$ is a viscosity subsolution of the HJB equation
(\ref{HJB}).
\end{proposition}

\textit{Proof. }Arguing by contradiction, we assume that $V$ is not a
subsolution of (\ref{HJB}) at $(x_{0},y_{0})$ with $x_{0}>0$ and $y_{0}>0$.
With a similar proof to the one of Proposition 3.1 of Azcue and Muler (2014),
but extending the definitions to two variables, we first show that there exist
$\varepsilon>0$, $h\in(0,\min\{x_{0}/2,y_{0}/2\})$ and a continuously
differentiable function $\psi:\mathbf{R}_{+}^{2}\rightarrow\mathbf{R}$ such
that $\psi$ is a test function for the subsolution of Equation (\ref{HJB}) at
$\left(  x_{0},y_{0}\right)  $ and satisfies%

\begin{equation}
\psi_{x}(x,y)\geq a_{1},\psi_{y}(x,y)\geq a_{2} \label{Desig1aprima}%
\end{equation}
for $(x,y)\in\lbrack0,x_{0}+h]\times\lbrack0,y_{0}+h],$
\begin{equation}
\mathcal{L}\left(  \psi\right)  (x,y)\leq-2\varepsilon\delta\label{Desig1a}%
\end{equation}
for $(x,y)\in\lbrack x_{0}-h,x_{0}+h]\times\lbrack y_{0}-h,y_{0}+h]$, and%

\begin{equation}
V(x,y)\leq\psi(x,y)-2\varepsilon\label{Desig2a}%
\end{equation}
for $\left(  x,y\right)  \in\mathbf{R}_{+}^{2}\setminus(x_{0}-h/2,x_{0}%
+h/2)\times(y_{0}-h/2,y_{0}+h/2)$.

Since $\psi$ is continuously differentiable, we can find a positive constant
$C$ such that
\begin{equation}
\mathcal{L}(\psi)(x,y)\leq C \label{Desig3a}%
\end{equation}
for all $(x,y)\in\lbrack0,x_{0}+2h]\times\lbrack0,y_{0}+2h]$.

Consider
\[
0<\theta<\left\{  \frac{h}{2\max\left\{  p_{1},p_{2}\right\}  }\text{, }%
\frac{\lambda}{4\delta\left(  \delta+\lambda\right)  },\frac{\varepsilon
\lambda}{2C\left(  \delta+\lambda\right)  }\right\}  \text{,}%
\]
and let us take any admissible strategy $\overline{L}\in\Pi_{x_{0},y_{0}}$.
Consider the corresponding controlled risk process $\left(  X_{t}%
,Y_{t}\right)  $ starting at $(x_{0},y_{0})$, and define the stopping times%
\begin{align*}
\tau^{b}  &  =\inf\{t>0:\text{ }\left(  X_{t},Y_{t}\right)  \in\partial\left(
\lbrack x_{0}-h,x_{0}+h]\times\lbrack y_{0}-h,y_{0}+h]\right)  \}\text{, }\\
\underline{\tau}  &  =\inf\{t>0:\text{ }\left(  X_{t},Y_{t}\right)
\in\mathbf{R}_{+}^{2}-[x_{0}-h,x_{0}+h]\times\lbrack y_{0}-h,y_{0}+h]\}
\end{align*}
and $\tau^{\ast}=\tau^{b}\wedge(\underline{\tau}+\theta)\wedge\overline{\tau}%
$. Note that $\tau^{\ast}$ is finite for $h$ small enough and that it is
necessary to introduce $\theta$ because before a lump sum dividend payment,
$(X_{\underline{\tau}},Y_{\underline{\tau}})$ can be in $[x_{0}-h,x_{0}%
+h]\times\lbrack y_{0}-h,y_{0}+h]$ and $(X_{\underline{\tau}},Y_{\underline
{\tau}})\in\mathbf{R}_{+}^{2}-[x_{0}-h,x_{0}+h]\times\lbrack y_{0}%
-h,y_{0}+h].$

Let us show that%

\begin{equation}
V(X_{\tau^{\ast}},Y_{\tau^{\ast}})\leq\psi(X_{\tau^{\ast}},Y_{\tau^{\ast}%
})-2\varepsilon\label{Vborde}%
\end{equation}
if $\tau^{\ast}=\tau^{b}\wedge(\underline{\tau}+\theta)<\overline{\tau}$.
There are two possibilities:

(1) If $\tau^{\ast}=\tau^{b},$ $\left(  X_{\tau^{\ast}},Y_{\tau^{\ast}%
}\right)  \in\partial\left(  \lbrack x_{0}-h,x_{0}+h]\times\lbrack
y_{0}-h,y_{0}+h]\right)  $ and so, from (\ref{Desig2a}), we obtain
$V(X_{\tau^{\ast}},Y_{\tau^{\ast}})\leq\psi(X_{\tau^{\ast}},Y_{\tau^{\ast}%
})-2\varepsilon$,

(2) If $\tau^{\ast}=\underline{\tau}+\theta,$ the distance from $(X_{\tau
^{\ast}},Y_{\tau^{\ast}})$ to $(x_{0},y_{0})$ is at least $h/2\geq
h-\max\left\{  p_{1},p_{2}\right\}  \theta$, so from (\ref{Desig2a}), we get
(\ref{Vborde}).

Note that $(X_{s^{-}},Y_{s^{-}})\in\lbrack0,x_{0}+h+p_{1}\theta]\times
\lbrack0,y_{0}+h+p_{2}\theta)]$ $\subset\lbrack0,x_{0}+2h]\times\lbrack
0,y_{0}+2h)]$ for $s\leq\tau^{\ast}$, so we have that%
\[
\mathcal{L}(\psi)(X_{s^{-}},Y_{s^{-}})\leq C\text{ for }s\leq\tau^{\ast
}\text{.}%
\]
Since $L_{t}^{i}$ , with $i=1,2$,$\ $is non-decreasing and left continuous, it
can be written as
\begin{equation}
L_{t}^{i}=\int\nolimits_{0}^{t}dL_{s}^{i,c}+\sum_{\substack{X_{s^{+}}\neq
X_{s} \\s<t}}(L_{s^{+}}^{i}-L_{s}^{i}) , \label{Lt}%
\end{equation}
where $L_{s}^{i,c}$ is a continuous and non-decreasing function. Since the
function $\psi$ is continuously differentiable in $\mathbf{R}_{+}^{2}$, using
the expression (\ref{Lt}) and the change of variables formula for finite
variation processes, we can write%
\[%
\begin{array}
[c]{l}%
\psi(X_{\tau^{\ast}},Y_{\tau^{\ast}})e^{-\delta\tau^{\ast}}-\psi(x_{0}%
,y_{0})\\%
\begin{array}
[c]{ll}%
= & \int\nolimits_{0}^{\tau^{\ast}}\left(  p_{1}\psi_{x}(X_{s^{-}},Y_{s^{-}%
})+p_{2}\psi_{y}(X_{s^{-}},Y_{s^{-}})\right)  e^{-\delta s}ds\\
& +\sum\limits_{\substack{X_{s^{-}}\neq X_{s} \\s\leq\tau^{\ast}}}\left(
\psi(X_{s},Y_{s})-\psi(X_{s^{-}},Y_{s^{-}})\right)  e^{-\delta s}%
+\sum\limits_{\substack{Y_{s^{-}}\neq Y_{s} \\s\leq\tau^{\ast}}}\left(
\psi(X_{s},Y_{s})-\psi(X_{s^{-}},Y_{s^{-}})\right)  e^{-\delta s}\\
& -\int\nolimits_{0}^{\tau^{\ast}}\psi_{x}(X_{s^{-}},Y_{s^{-}})e^{-\delta
s}dL_{s}^{1,c}+\sum\limits_{\substack{X_{s^{+}}\neq X_{s} \\s<\tau^{\ast}%
}}\left(  \psi(X_{s^{+}},Y_{s^{+}})-\psi(X_{s},Y_{s})\right)  e^{-\delta s}\\
& -\int\nolimits_{0}^{\tau^{\ast}}\psi_{y}(X_{s^{-}},Y_{s^{-}})e^{-\delta
s}dL_{s}^{2,c}+\sum\limits_{\substack{Y_{s^{+}}\neq Y_{s} \\s<\tau^{\ast}%
}}\left(  \psi(X_{s^{+}},Y_{s^{+}})-\psi(X_{s},Y_{s})\right)  e^{-\delta s}\\
& -\delta\int\nolimits_{0}^{\tau^{\ast}}\psi(X_{s^{-}},Y_{s^{-}})e^{-\delta
s}ds.
\end{array}
\end{array}
\]

Note that $(X_{s},Y_{s})\in\mathbf{R}_{+}^{2}$ for $s\leq\tau^{\ast}$ except
in the case that $\tau^{\ast}=\overline{\tau},$ where $X_{\tau^{\ast}}%
+Y_{\tau^{\ast}}<0$. Here we are extending the definition of $\psi$ as
\[
\psi(x,y)=a_{1}V_{1}^{0}(x)I_{x\geq0}+a_{2}V_{2}^{0}(y)I_{y\geq0}%
\]
for $x+y<0$. We have that $X_{s^{+}}\neq X_{s}$ only at the jumps of
$L_{s}^{1}$, and in this case $X_{s^{+}}-X_{s}=-\left(  L_{s^{+}}^{1}%
-L_{s}^{1}\right)  $. Since $\overline{L}$ is admissible we have that
$X_{s^{+}}=X_{s}-\left(  L_{s^{+}}^{1}-L_{s}^{1}\right)  \geq0$. We can write%

\begin{equation}%
\begin{array}
[c]{l}%
-\int\nolimits_{0}^{\tau^{\ast}}\psi_{x}(X_{s^{-}},Y_{s^{-}})e^{-\delta
s}dL_{s}^{1,c}+\sum\limits_{\substack{X_{s^{+}}\neq X_{s}\\s<\tau^{\ast}%
}}\left(  \psi(X_{s^{+}},Y_{s^{+}})-\psi(X_{s},Y_{s})\right)  e^{-\delta s}\\%
\begin{array}
[c]{cl}%
= & -\int\nolimits_{0}^{\tau^{\ast}}\psi_{x}(X_{s^{-}},Y_{s^{-}})e^{-\delta
s}dL_{s}^{1,c}+\sum\limits_{\substack{X_{s^{+}}\neq X_{s}\\s<\tau^{\ast}%
}}\left(  \int\nolimits_{0}^{L_{s^{+}}^{1}-L_{s}^{1}}\psi_{x}(X_{s}%
-\alpha,Y_{s})d\alpha\right)  e^{-\delta s}\\
\leq & -\int\nolimits_{0}^{\tau^{\ast}}a_{1}e^{-\delta s}dL_{s}^{1,c}%
-a_{1}\sum\limits_{\substack{L_{s^{+}}^{1}\neq L_{s}^{1}\\s<\tau^{\ast}%
}}\left(  \int\nolimits_{0}^{L_{s^{+}}^{1}-L_{s}^{1}}d\alpha\right)
e^{-\delta s}\\
= & -a_{1}\int_{0}^{\tau^{\ast}}e^{-\delta s}dL_{s}^{1}.
\end{array}
\end{array}
\label{Paso 2}%
\end{equation}

Similarly,%

\begin{equation}%
\begin{array}
[c]{l}%
-\int\nolimits_{0}^{\tau^{\ast}}\psi_{x}(X_{s^{-}},Y_{s^{-}})e^{-\delta
s}dL_{s}^{2,c}+\sum\limits_{\substack{X_{s^{+}}\neq X_{s}\\s<\tau^{\ast}%
}}\left(  \psi(X_{s^{+}},Y_{s^{+}})-\psi(X_{s},Y_{s})\right)  e^{-\delta s}\\%
\begin{array}
[c]{cl}%
= & -\int\nolimits_{0}^{\tau^{\ast}}\psi_{x}(X_{s^{-}},Y_{s^{-}})e^{-\delta
s}dL_{s}^{2,c}-\sum\limits_{\substack{L_{s^{+}}^{2}\neq L_{s}^{2}%
\\s<\tau^{\ast}}}\left(  \int\nolimits_{0}^{L_{s^{+}}^{2}-L_{s}^{2}}\psi
_{x}(X_{s},Y_{s}-\alpha)d\alpha\right)  e^{-\delta s}\\
\leq & -a_{2}\int_{0}^{\tau^{\ast}}e^{-\delta s}dL_{s}^{2}.
\end{array}
\end{array}
\label{Paso2Y}%
\end{equation}
On the other hand, $X_{s}\neq X_{s^{-}}$ only at the arrival of a claim for
Company One, so%

\begin{equation}%
\begin{array}
[c]{cll}%
M_{t}^{1} & = &
%TCIMACRO{\tsum \limits_{\substack{X_{s^{-}}\neq X_{s}\\s\leq t}}}%
%BeginExpansion
{\textstyle\sum\limits_{\substack{X_{s^{-}}\neq X_{s}\\s\leq t}}}
%EndExpansion
\left(  \psi(X_{s},Y_{s})-\psi(X_{s^{-}},Y_{s^{-}})\right)  e^{-\delta s}\\
&  & -\lambda_{1}\int\nolimits_{0}^{t}e^{-\delta s}\int\nolimits_{0}%
^{X_{s^{-}}}\left(  \psi(X_{s^{-}}-\alpha,Y_{s^{-}})-\psi(X_{s^{-}},Y_{s^{-}%
})\right)  dF^{1}(\alpha)ds\\
&  & -\lambda_{1}\int\nolimits_{0}^{t}e^{-\delta s}\int\nolimits_{X_{s^{-}}%
}^{X_{s^{-}}+Y_{s^{-}}}\left(  \psi(0,X_{s^{-}}+Y_{s^{-}}-\alpha
)-\psi(X_{s^{-}},Y_{s^{-}})\right)  dF^{1}(\alpha)ds\\
&  & -\lambda_{1}\int\nolimits_{0}^{t}e^{-\delta s}\int\nolimits_{X_{s^{-}%
}+Y_{s^{-}}}^{\infty}\left(  a_{2}V_{2}^{0}(Y_{s^{-}})-\psi(X_{s^{-}}%
,Y_{s^{-}})\right)  dF^{1}(\alpha)ds
\end{array}
\label{M1t}%
\end{equation}
is a martingale with zero expectation for $t\leq\overline{\tau}$. Analogously,%

\begin{equation}%
\begin{array}
[c]{cll}%
M_{t}^{2} & = &
%TCIMACRO{\tsum \limits_{\substack{Y_{s^{-}}\neq Y_{s}\\s\leq t}}}%
%BeginExpansion
{\textstyle\sum\limits_{\substack{Y_{s^{-}}\neq Y_{s}\\s\leq t}}}
%EndExpansion
\left(  \psi(X_{s},Y_{s})-\psi(X_{s^{-}},Y_{s^{-}})\right)  e^{-\delta s}\\
&  & -\lambda_{2}\int\nolimits_{0}^{t}e^{-\delta s}\int\nolimits_{0}%
^{Y_{s^{-}}}\left(  \psi(X_{s^{-}},Y_{s^{-}}-\alpha)-\psi(X_{s^{-}},Y_{s^{-}%
})\right)  dF^{2}(\alpha)ds\\
&  & -\lambda_{2}\int\nolimits_{0}^{t}e^{-\delta s}\int\nolimits_{X_{s^{-}}%
}^{X_{s^{-}}+Y_{s^{-}}}\left(  \psi(X_{s^{-}}+Y_{s^{-}}-\alpha,0)-\psi
(X_{s^{-}},Y_{s^{-}})\right)  dF^{2}(\alpha)ds\\
&  & -\lambda_{2}\int\nolimits_{0}^{t}e^{-\delta s}\int\nolimits_{X_{s^{-}%
}+Y_{s^{-}}}^{\infty}\left(  a_{1}V_{1}^{0}(Y_{s^{-}})-\psi(X_{s^{-}}%
,Y_{s^{-}})\right)  dF^{2}(\alpha)ds
\end{array}
\label{M2t}%
\end{equation}
is also a martingale with zero expectation for $t\leq\overline{\tau}$. So we get%

\begin{equation}%
\begin{array}
[c]{ccc}%
\psi(X_{\tau^{\ast}},Y_{\tau^{\ast}})e^{-\delta\tau^{\ast}}-\psi(x_{0}%
,y_{0}) & \leq & \int\nolimits_{0}^{\tau^{\ast}}\mathcal{L}(\psi)(X_{s^{-}%
},Y_{s^{-}})e^{-\delta s}+M_{\tau^{\ast}}^{1}+M_{\tau^{\ast}}^{2}\\
&  & -a_{1}\int_{0}^{\tau^{\ast}}e^{-\delta s}dL_{s}^{1}-a_{2}\int_{0}%
^{\tau^{\ast}}e^{-\delta s}dL_{s}^{2}.
\end{array}
\label{nueva1}%
\end{equation}
Using the second inequality of (\ref{Desig1a}), (\ref{Desig3a}) and the
definition of $\theta$ we get
\begin{equation}%
\begin{array}
[c]{ccl}%
\int\nolimits_{0}^{\tau^{\ast}}\mathcal{L}(\psi)(X_{s^{-}},Y_{s^{-}%
})e^{-\delta s}ds & \leq & \int\nolimits_{0}^{\tau^{b}\wedge\underline{\tau
}\wedge\overline{\tau}}\mathcal{L}(\psi)(X_{s^{-}},Y_{s^{-}})e^{-\delta
s}ds+C\theta\\
& \leq & -2\varepsilon\delta\int\nolimits_{0}^{\tau^{b}\wedge\underline{\tau
}\wedge\overline{\tau}}e^{-\delta s}ds+C\theta\\
& \leq & -2\varepsilon\delta\int\nolimits_{0}^{\tau^{\ast}}e^{-\delta
s}ds+I_{\tau^{b}\wedge\underline{\tau}\wedge\overline{\tau}<\tau^{\ast}%
}2\varepsilon\delta\int\nolimits_{\tau^{b}\wedge\underline{\tau}%
\wedge\overline{\tau}}^{\tau^{\ast}}e^{-\delta s}ds+C\theta\\
& \leq & -2\varepsilon(1-e^{-\delta\tau^{\ast}})+\varepsilon\lambda/\left(
\delta+\lambda\right)  .
\end{array}
\label{nueva2}%
\end{equation}
From (\ref{Desig}), Lemma \ref{DPP}, (\ref{Vborde}), (\ref{nueva1}) and
(\ref{nueva2}), it follows that%

\begin{equation}%
\begin{array}
[c]{l}%
V(x_{0},y_{0})\\%
\begin{array}
[c]{ll}%
= & \sup\nolimits_{\overline{L}}(E_{x_{0,}y_{0}}(a_{1}\int_{0}^{\tau^{\ast}%
}e^{-\delta s}dL_{s}^{1}+a_{2}\int_{0}^{\tau^{\ast}}e^{-\delta s}dL_{s}%
^{2}+e^{-\delta\tau^{\ast}}V(X_{\tau^{\ast}},Y_{\tau^{\ast}})I_{\tau^{\ast
}<\overline{\tau}}\\
& +\left(  a_{1}V_{1}^{0}(X_{\overline{\tau}}^{\overline{L}})+a_{2}V_{2}%
^{0}(Y_{\overline{\tau}}^{\overline{L}})\right)  e^{-\delta\tau^{\ast}}%
I_{\tau^{\ast}=\overline{\tau}}))\\
\leq & \sup\nolimits_{\overline{L}}(E_{x_{0,}y_{0}}(a_{1}\int_{0}^{\tau^{\ast
}}e^{-\delta s}dL_{s}^{1}+a_{2}\int_{0}^{\tau^{\ast}}e^{-\delta s}dL_{s}%
^{2}+e^{-\delta\tau^{\ast}}\left(  \psi(X_{\tau^{\ast}},Y_{\tau^{\ast}%
})-2\varepsilon\right)  I_{\tau^{\ast}<\overline{\tau}}\\
& +\left(  a_{1}V_{1}^{0}(X_{\overline{\tau}}^{\overline{L}})+a_{2}V_{2}%
^{0}(Y_{\overline{\tau}}^{\overline{L}})\right)  e^{-\delta\tau^{\ast}}%
I_{\tau^{\ast}=\overline{\tau}}))\\
\leq & \sup\nolimits_{\overline{L}}E_{x_{0},y_{0}}\left(  \int\nolimits_{0}%
^{\tau^{\ast}}\mathcal{L}\left(  \psi\right)  (X_{s^{-}},Y_{s^{-}})e^{-\delta
s}ds+M_{\tau^{\ast}}^{1}+M_{\tau^{\ast}}^{2}-2\varepsilon e^{-\delta\tau
^{\ast}}I_{\tau^{\ast}<\overline{\tau}}+\psi(x_{0},y_{0})\right) \\
\leq & \psi(x_{0},y_{0})-\varepsilon\lambda/\left(  \delta+\lambda\right) \\
< & \psi(x_{0},y_{0})
\end{array}
\end{array}
\label{Desig}%
\end{equation}
and this contradicts the assumption that $V(x_{0},y_{0})=\psi(x_{0},y_{0}%
)$.\hfill{\small $\square$}

From the above two propositions we get the following result.

\begin{corollary}
$V$ is a viscosity solution of of the HJB equation (\ref{HJB}).
\end{corollary}

\section{Smallest Viscosity Solution}

\label{sec4}

Let us prove now that the optimal value function $V$ is the smallest viscosity
supersolution of (\ref{HJB}).

We say that the function $u:\mathbf{R}_{+}^{2}\rightarrow\mathbf{R}$ satisfies
the growth condition A.1, if%

\[
u(x,y)\leq K+a_{1}x+a_{2}y\text{ for all }(x,y)\in\mathbf{R}_{+}^{2}\text{.}%
\]

The following Lemma is technical and will be used to prove Proposition
\ref{MenorSuper}.

\begin{lemma}
\label{A.1}Fix $x_{0}>0$ and $y_{0}>0$ and let $\overline{u}$\ be a
non-negative supersolution of (\ref{HJB}) satisfying the growth condition
A.1.\ We can find a sequence of positive functions $\overline{u}%
_{m}:\mathbf{R}_{+}^{2}\rightarrow\mathbf{R}$\ such that:

(a) $\overline{u}_{m}$\ is continuously differentiable.

(b) $\overline{u}_{m}\ $satisfies the growth condition A.1.

(c)$\ p\leq p_{1}\overline{u}_{m,x}+p_{2}\overline{u}_{m,y}$\ $\leq\left(
\delta+\lambda\right)  \overline{u}_{m}$ in $\mathbf{R}_{+}^{2}$.

(d) $\overline{u}_{m}$\ $\searrow$ $\overline{u}$\ uniformly on compact sets
in $\mathbf{R}_{+}^{2}$ and $\nabla\overline{u}_{m}$\ converges to
$\nabla\overline{u}$\ a.e. in $\mathbf{R}_{+}^{2}$.

(e) There exists a sequence $c_{m}$ with $\lim\limits_{m\rightarrow\infty
}c_{m}=0$ such that
\[
\sup\nolimits_{\left(  x,y\right)  \in A_{0}}\mathcal{L}(\overline{u}%
_{m})\left(  x,y\right)  \leq c_{m},\text{where }A_{0}=[0,x_{0}]\times
\lbrack0,y_{0}].
\]

\end{lemma}

\textit{Proof.} The proof follows by standard convolution arguments and is the
extension to two variables of Lemma 4.1 in Azcue an\ Muler \cite{azmu}%
.\hfill{\small $\square$ }

\begin{proposition}
\smallskip\label{MenorSuper} The optimal value function $V$ is the smallest
viscosity supersolution of (\ref{HJB}) satisfying growth condition A.1.
\end{proposition}

\textit{Proof. }Let $\overline{u}$\ be a non-negative supersolution of
(\ref{HJB}) satisfying the growth condition A.1 and let $\overline{L}\in
\Pi_{x,y}$; define $\left(  X_{t},Y_{t}\right)  $ as the corresponding
controlled risk process starting at $\left(  x,y\right)  $. Consider the
function $\overline{u}_{m}$ of Lemma \ref{A.1} in $\mathbf{R}_{+}^{2}$ ; we
extend this function as
\[
\overline{u}_{m}(x,y)=a_{1}V_{1}^{0}(x)I_{x\geq0}+a_{2}V_{2}^{0}(y)I_{y\geq
0}\text{ for }x+y<0.
\]
As in the proof of Proposition \ref{Prop V Subsolution}, we get%

\begin{equation}%
\begin{array}
[c]{l}%
\overline{u}_{m}(X_{t\wedge\overline{\tau}},Y_{t\wedge\overline{\tau}%
})e^{-\delta(t\wedge\overline{\tau})}-\overline{u}_{m}(x,y)\\
\leq\int\nolimits_{0}^{t\wedge\overline{\tau}}\mathcal{L}(\overline{u}%
_{m})(X_{s^{-}},Y_{s^{-}})e^{-\delta s}ds-a_{1}\int\nolimits_{0}%
^{t\wedge\overline{\tau}}e^{-\delta s}dL_{s}^{1}-a_{2}\int\nolimits_{0}%
^{t\wedge\overline{\tau}}e^{-\delta s}dL_{s}^{2}+M_{t\wedge\overline{\tau}%
}^{1}+M_{t\wedge\overline{\tau}}^{2}\text{,}%
\end{array}
\label{ItoUnMenorSuper}%
\end{equation}
where $M_{t}^{1}\ $and $M_{t}^{2}$ are zero-expectation martingales. So we
obtain that%

\[%
\begin{array}
[c]{l}%
\overline{u}_{m}(X_{t\wedge\overline{\tau}},Y_{t\wedge\overline{\tau}%
})e^{-\delta(t\wedge\overline{\tau})}I_{\overline{\tau}>t}-\overline{u}%
_{m}(x,y)\\%
\begin{array}
[c]{ll}%
\leq & \int\nolimits_{0}^{t\wedge\overline{\tau}}\mathcal{L}(\overline{u}%
_{m})(X_{s^{-}},Y_{s^{-}})e^{-\delta s}ds+M_{t\wedge\overline{\tau}}%
^{1}+M_{t\wedge\overline{\tau}}^{2}\\
&
\begin{array}
[c]{l}%
-a_{1}\left(  \int_{0}^{t\wedge\overline{\tau}}e^{-\delta s}dL_{s}%
^{1}+e^{-\delta\left(  t\wedge\overline{\tau}\right)  }V_{1}^{0}%
(X_{t\wedge\overline{\tau}})I_{\overline{\tau}\leq t}\right) \\
-a_{2}\left(  \int_{0}^{t\wedge\overline{\tau}}e^{-\delta s}dL_{s}%
^{2}+e^{-\delta\left(  t\wedge\overline{\tau}\right)  }V_{2}^{0}%
(Y_{t\wedge\overline{\tau}})I_{\overline{\tau}\leq t}\right)  .
\end{array}
\end{array}
\end{array}
\]

Using that$\mathcal{\ }$both $L_{t}^{1}$ and $L_{t}^{2}$ are non-decreasing
processes, from the monotone convergence theorem we get
\begin{equation}%
\begin{array}
[c]{l}%
\lim\limits_{t\rightarrow\infty}(E_{x,y}(a_{1}\left(  \int_{0}^{t\wedge
\overline{\tau}}e^{-\delta s}dL_{s}^{1}+e^{-\delta\left(  t\wedge
\overline{\tau}\right)  }V_{1}^{0}(X_{t\wedge\overline{\tau}})I_{\overline
{\tau}\leq t}\right) \\
+a_{2}\left(  \int_{0}^{t\wedge\overline{\tau}}e^{-\delta s}dL_{s}%
^{2}+e^{-\delta\left(  t\wedge\overline{\tau}\right)  }V_{2}^{0}%
(Y_{t\wedge\overline{\tau}})I_{\overline{\tau}\leq t}\right)  ))\\
=V_{\overline{L}}(x,y).
\end{array}
\label{monotone1}%
\end{equation}
From Lemma \ref{A.1}(c), we have%

\[
-(\delta+\lambda)\overline{u}_{m}(x,y)\leq\mathcal{L}(\overline{u}%
_{m})(x,y)\leq\lambda\overline{u}_{m}(x,y)+U(x,y).
\]
But using Lemma \ref{A.1}(b) and the inequality $X_{s}\leq x+p_{1}s$,
$Y_{s}\leq y+p_{2}s$ we get
\begin{equation}
\overline{u}_{m}(X_{s},Y_{s})\leq K+a_{1}X_{s}+a_{2}Y_{s}\leq K+a_{1}%
x+a_{2}y+ps. \label{Acotacionu}%
\end{equation}
So, using the bounded convergence theorem, we obtain
\begin{equation}
\lim\limits_{t\rightarrow\infty}E_{x,y}\left(  \int\nolimits_{0}%
^{t\wedge\overline{\tau}}\mathcal{L}(\overline{u}_{m})(X_{s^{-}},Y_{s^{-}%
})e^{-\delta s}ds\right)  =E_{x,y}\left(  \int\nolimits_{0}^{\overline{\tau}%
}\mathcal{L}(\overline{u}_{m})(X_{s^{-}},Y_{s^{-}})e^{-\delta s}ds\right)  .
\label{monotone2}%
\end{equation}
From (\ref{ItoUnMenorSuper}), (\ref{monotone1}) and (\ref{monotone2}), we get
\begin{equation}
\lim\limits_{t\rightarrow\infty}E_{x,y}\left(  \overline{u}_{m}(X_{t\wedge
\overline{\tau}},Y_{t\wedge\overline{\tau}})e^{-\delta(t\wedge\overline{\tau
})}I_{\overline{\tau}<t}\right)  -\overline{u}_{m}(x,y)\leq E_{x,y}\left(
\int\nolimits_{0}^{\overline{\tau}}\mathcal{L}(\overline{u}_{m})(X_{s^{-}%
},Y_{s^{-}})e^{-\delta s}ds\right)  -V_{\overline{L}}(x,y). \label{limite0}%
\end{equation}

Next, we show that
\begin{equation}
\lim\limits_{t\rightarrow\infty}E_{x,y}\left(  \overline{u}_{m}(X_{t\wedge
\overline{\tau}},Y_{t\wedge\overline{\tau}})e^{-\delta(t\wedge\overline{\tau
})}I_{\overline{\tau}>t}\right)  =0. \label{limite1}%
\end{equation}
From (\ref{Acotacionu}), there exists a $\overline{K}$ such that
\[%
\begin{array}
[c]{ccl}%
E_{x,y}\left(  \overline{u}_{m}(X_{t\wedge\overline{\tau}},Y_{t\wedge
\overline{\tau}})e^{-\delta(t\wedge\overline{\tau})}I_{\overline{\tau}%
>t}\right)  & \leq & \left(  \overline{K}+a_{1}x+a_{2}y+pt\right)  e^{-\delta
t}\text{.}%
\end{array}
\]
Since the last expression goes to $0$ as $t$ goes to infinity, we have
(\ref{limite1}). Let us prove now that
\begin{equation}
\limsup\limits_{m\rightarrow\infty}E_{x,y}\left(  \int\nolimits_{0}%
^{\overline{\tau}}\mathcal{L}(\overline{u}_{m})(X_{s^{-}},Y_{s^{-}})e^{-\delta
s}ds\right)  \leq0. \label{limite2}%
\end{equation}
Given any $\varepsilon>0,$ we can find $T$ such that
\begin{equation}
\int\nolimits_{T}^{\infty}\mathcal{L}(\overline{u}_{m})(X_{s^{-}},Y_{s^{-}%
})e^{-\delta s}ds<\frac{\varepsilon}{2} \label{arreglado5}%
\end{equation}
for any $m\geq1$, as by virtue of (\ref{Acotacionu}), growth condition A.1,
Lemma \ref{A.1}(b) and Lemma \ref{A.1}(c), and the growth property of
$V_{1}^{0}$ and $V_{2}^{0},$ there exist positive constants $k_{0},k_{1}%
$,$k_{2}$ and $\overline{p}$ such that
\[%
\begin{array}
[c]{lll}%
\mathcal{L}(\overline{u}_{m})(X_{s^{-}},Y_{s^{-}}) & \leq & \lambda
\overline{u}_{m}(X_{s^{-}},Y_{s^{-}})+U(X_{s^{-}},Y_{s^{-}})\\
& \leq & k_{0}+k_{1}x+k_{2}y+\overline{p}s.
\end{array}
\]

Note that for $s\leq T$, $X_{s^{-}}$ $\leq x_{0}:=x+p_{1}T$ , $Y_{s^{-}}$
$\leq y_{0}:=y+p_{2}T$. From Lemma \ref{A.1}(e) we can find $m_{0}$ large
enough such that for any $m\geq m_{0}$%
\[
\int\nolimits_{0}^{T}\mathcal{L}(\overline{u}_{m})(X_{s^{-}},Y_{s^{-}%
})e^{-\delta s}ds\leq c_{m}\int\nolimits_{0}^{T}e^{-\delta s}ds\leq\frac
{c_{m}}{\delta}\leq\frac{\varepsilon}{2}%
\]
and so we get (\ref{limite2}). Then, from (\ref{limite0}) and using
(\ref{limite1}) and (\ref{limite2}), we obtain%

\begin{equation}
\overline{u}(x,y)=\lim\limits_{m\rightarrow\infty}\overline{u}_{m}(x,y)\geq
V_{\overline{L}}(x,y). \label{supersolucionMayorqueVL}%
\end{equation}
Since $V$ is a viscosity solution of (\ref{HJB}), the result follows.\hfill
{\small $\square$}

From the previous proposition we can deduce the usual viscosity verification result.

\begin{corollary}
\label{VerificacionDividendos} Consider a family of admissible strategies
$\{\overline{L}^{x,y}\in\Pi_{x,y}:\left(  x,y\right)  \in\mathbf{R}_{+}^{2}%
\}$. If the function $V_{\overline{L}^{x,y}}(x,y)$ is a viscosity
supersolution of (\ref{HJB}) for all $\left(  x,y\right)  \in\mathbf{R}%
_{+}^{2}$, then $V_{\overline{L}^{x,y}}(x,y)$ is the optimal value function
(\ref{DefinitionV}).
\end{corollary}

\section{Iterative Approach\label{Iterative}}

\label{sec5}

In this section, we approximate the optimal value function $V$ defined in
(\ref{DefinitionV}) by an increasing sequence of value functions of strategies
which pay dividends (and collaborate if it is necessary) up to the $n$-th
claim (regardless from which company) and then follow the
take-the-money-and-run strategy. Given initial surplus levels $(x,y)$, the
take-the-money-and-run admissible strategy $\overline{L}^{0}$pays immediately
the entire surplus $x$ and $y$ as dividends (that is $X_{0^{+}}=Y_{0^{+}}=0$),
and then pays the incoming premium as dividends until the first claim, where
the company facing that claim gets ruined. Note that under this strategy the
companies can not help each other.

Consider $\tau_{n}$ as the time of arrival of the $n$-th claim regardless from
which company, that is the $n$-th point of the Poisson process $N_{t}%
=N_{t}^{1}+N_{t}^{2}$. We define the set $\Pi_{x,y}^{n}$ of all the admissible
strategies in $\Pi_{x,y}$ which follow $\overline{L}^{0}$ right after
$\tau_{n}$. Let us define%

\begin{equation}
V^{n\text{ }}(x,y)=\sup_{\overline{L}\in\Pi_{x,y}^{n}}V_{\overline{L}}(x,y)
\label{Definicion Vn}%
\end{equation}
for $n\geq1$, we also define $V^{0}=V_{\overline{L}^{0}}$. We can write%

\[%
\begin{array}
[c]{ccc}%
V^{0}(x,y) & = & a_{1}x+a_{2}y+\frac{\lambda_{1}}{\delta+\lambda}\left(
\frac{p}{\lambda}+a_{1}V_{2}^{0}(0)\right) \\
&  & +\frac{\lambda_{2}}{\delta+\lambda}\left(  \frac{p}{\lambda}+a_{2}%
V_{1}^{0}(0)\right)  .
\end{array}
\]
Note that, for $n\geq1$, we have%

\begin{equation}
V^{n\text{ }}(x,y)=\sup_{\overline{L}\in\Pi_{x,y}}V_{\overline{L}}^{n}(x,y),
\label{Step_n_Problem}%
\end{equation}
where
\begin{equation}%
\begin{array}
[c]{ll}%
V_{\overline{L}}^{n}(x,y)= & E_{x,y}(a_{1}\int_{0}^{\tau_{1}}e^{-\delta
s}dL_{s}^{1}+a_{2}\int_{0}^{\tau_{1}}e^{-\delta s}dL_{s}^{2}+e^{-\delta
\tau_{1}}V^{n-1}(X_{\tau_{1}}^{\overline{L}},Y_{\tau_{1}}^{\overline{L}%
})I_{\tau_{1}<\overline{\tau}}\\
& +e^{-\delta\overline{\tau}}\left(  a_{1}V_{1}^{0}(X_{\overline{\tau}%
}^{\overline{L}})+a_{2}V_{2}^{0}(Y_{\overline{\tau}}^{\overline{L}})\right)
I_{\tau_{1}=\overline{\tau}}).
\end{array}
\label{VnLbarra}%
\end{equation}

In this expression, we only consider the admissible strategy $\overline{L}%
\in\Pi_{x,y}$ for $t\leq\tau_{1}$. The following DPP holds.

\begin{lemma}
\label{DPPIterativa} For any initial surplus $\left(  x,y\right)  $ in
$\mathbf{R}_{+}^{2}$ and any stopping time $\tau$, we can write%
\[%
\begin{array}
[c]{l}%
\begin{array}
[c]{ll}%
V^{n}(x,y)= & \sup_{\overline{L}\in\Pi_{x,y}}(E_{x,y}(a_{1}\int_{0}^{\tau
_{1}\wedge\tau}e^{-\delta s}dL_{s}^{1}+a_{2}\int_{0}^{\tau_{1}\wedge\tau
}e^{-\delta s}dL_{s}^{2}\\
& +e^{-\delta\tau}V^{n}(X_{\tau}^{\overline{L}},Y_{\tau}^{\overline{L}%
})I_{\tau<\tau_{1}\wedge\overline{\tau}}+e^{-\delta\tau_{1}}V^{n-1}%
(X_{\tau_{1}}^{\overline{L}},Y_{\tau_{1}}^{\overline{L}})I_{\tau_{1}\wedge
\tau=\tau_{1}<\overline{\tau}}\\
& +e^{-\delta\overline{\tau}}\left(  a_{1}V_{1}^{0}(X_{\overline{\tau}%
}^{\overline{L}})+a_{2}V_{2}^{0}(Y_{\overline{\tau}}^{\overline{L}})\right)
I_{\tau_{1}\wedge\tau=\tau_{1}=\overline{\tau}})).
\end{array}
\end{array}
\]

\end{lemma}

\begin{proposition}
\label{Monotonicidad Vn} We have that $V^{1}\leq V^{2}\leq...\leq V$.
\end{proposition}

\textit{Proof.} We prove the result by induction:

(a) $V^{0}\leq V^{1}$, because the strategy $\overline{L}^{0}\in\Pi_{x,y}^{1}$.

(b) Assume that $V^{n-2\text{ }}\leq V^{n-1\text{ }}$. By
(\ref{Step_n_Problem}), we have%

\[%
\begin{array}
[c]{lll}%
V^{n\text{ }}(x,y) & \geq & \sup_{\overline{L}\in\Pi_{x,y}}E_{x,y}(a_{1}%
\int_{0}^{\tau_{1}}e^{-\delta s}dL_{s}^{1}+a_{2}\int_{0}^{\tau_{1}}e^{-\delta
s}dL_{s}^{2}+e^{-\delta\tau_{1}}V^{n-2\text{ }}(X_{\tau_{1}}^{\overline{L}%
},Y_{\tau_{1}}^{\overline{L}})I_{\tau_{1}<\overline{\tau}}\\
&  &
\begin{array}
[c]{l}%
+e^{-\delta\overline{\tau}}\left(  a_{1}V_{1}^{0}(X_{\overline{\tau}%
}^{\overline{L}})+a_{2}V_{2}^{0}(Y_{\overline{\tau}}^{\overline{L}})\right)
I_{\tau_{1}=\overline{\tau}})\\
.
\end{array}
\\
& = & V^{n-1\text{ }}(x,y)~\text{.}\hfill\smallskip\ {\small \square}%
\end{array}
\]

The HJB equation for $V^{n\text{ }}$is given by%

\begin{equation}
\max\left\{  \mathcal{L}^{n}(V^{n})(x,y),a_{1}-V_{x}^{n}(x,y),a_{2}-V_{y}%
^{n}(x,y)\right\}  =0, \label{HJB Iterativa}%
\end{equation}
where%

\begin{equation}%
\begin{array}
[c]{lll}%
\mathcal{L}^{n}(V^{n})(x,y) & = & p_{1}V_{x}^{n}(x,y)+p_{2}V_{y}%
^{n}(x,y)-\left(  \delta+\lambda\right)  V^{n}(x,y)+\mathcal{I}(V^{n-1}%
)(x,y)+U(x,y).
\end{array}
\label{DefLn}%
\end{equation}

The following basic results about regularity and growth at infinity of $V^{n}
$, $n\geq1$ are similar to those of Lemmas \ref{V_GrowthCondition} and

\ref{V increasing_LocallyLip}.

\begin{lemma}
\label{Vn increasing_LocallyLip}The optimal value function $V^{n}$ satisfies
growth condition A.1, it is increasing and locally Lipschitz in $\mathbf{R}%
_{+}^{2}$ with
\[
a_{1}h\leq V^{n}(x+h,y)-V^{n}(x,y)\leq(e^{(\delta+\lambda)h/p_{1}}%
-1)V^{n}(x,y)
\]%
\[
a_{2}h\leq V^{n}(x,y+h)-V^{n}(x,y)\leq(e^{(\delta+\lambda)h/p_{2}}%
-1)V^{n}(x,y)
\]
for any $h>0\ $and for any $(x,y)$ in $\mathbf{R}_{+}^{2}$.
\end{lemma}

In the next two propositions, we see that $V^{n}$ is a viscosity solution of
the corresponding HJB equation.

\begin{proposition}
$V^{n}$ is a viscosity supersolution of the HJB equation (\ref{HJB Iterativa})
for $x>0\ $and $y>0$.
\end{proposition}

\textit{Proof. }Similar to the one given in Proposition
\ref{Prop V is a viscosity supersolution}.\hfill${\small \square}$

\begin{proposition}
\label{Prop V SubsolutionIterativa}$V^{n}$ is a viscosity subsolution of the
corresponding HJB equation (\ref{HJB Iterativa}).
\end{proposition}

\textit{Proof. }The proof of this proposition is similar to the one of
Proposition\textit{\ }\ref{Prop V Subsolution}, but using as martingales with
zero expectation%

\begin{equation}%
\begin{array}
[c]{cll}%
M_{t}^{3} & = &
%TCIMACRO{\tsum \limits_{\substack{X_{s^{-}}\neq X_{s}\\s\leq t}}}%
%BeginExpansion
{\textstyle\sum\limits_{\substack{X_{s^{-}}\neq X_{s}\\s\leq t}}}
%EndExpansion
\left(  V^{n-1}(X_{s},Y_{s})-\psi(X_{s^{-}},Y_{s^{-}})\right)  e^{-\delta s}\\
&  & -\lambda_{1}\int\nolimits_{0}^{t}e^{-\delta s}\int\nolimits_{0}%
^{X_{s^{-}}}\left(  V^{n-1}(X_{s^{-}}-\alpha,Y_{s^{-}})-\psi(X_{s^{-}%
},Y_{s^{-}})\right)  dF^{1}(\alpha)ds\\
&  & -\lambda_{1}\int\nolimits_{0}^{t}e^{-\delta s}\int\nolimits_{X_{s^{-}}%
}^{X_{s^{-}}+Y_{s^{-}}}\left(  V^{n-1}(0,X_{s^{-}}+Y_{s^{-}}-\alpha
)-\psi(X_{s^{-}},Y_{s^{-}})\right)  dF^{1}(\alpha)ds\\
&  & -\lambda_{1}\int\nolimits_{0}^{t}e^{-\delta s}\int\nolimits_{X_{s^{-}%
}+Y_{s^{-}}}^{\infty}\left(  a_{2}V_{2}^{0}(Y_{s^{-}})-\psi(X_{s^{-}}%
,Y_{s^{-}})\right)  dF^{1}(\alpha)ds
\end{array}
\label{Martingala3}%
\end{equation}
and%

\begin{equation}%
\begin{array}
[c]{cll}%
M_{t}^{4} & = &
%TCIMACRO{\tsum \limits_{\substack{Y_{s^{-}}\neq Y_{s}\\s\leq t}}}%
%BeginExpansion
{\textstyle\sum\limits_{\substack{Y_{s^{-}}\neq Y_{s}\\s\leq t}}}
%EndExpansion
\left(  V^{n-1}(X_{s},Y_{s})-\psi(X_{s^{-}},Y_{s^{-}})\right)  e^{-\delta s}\\
&  & -\lambda_{2}\int\nolimits_{0}^{t}e^{-\delta s}\int\nolimits_{0}%
^{Y_{s^{-}}}\left(  V^{n-1}(X_{s^{-}},Y_{s^{-}}-\alpha)-\psi(X_{s^{-}%
},Y_{s^{-}})\right)  dF^{2}(\alpha)ds\\
&  & -\lambda_{2}\int\nolimits_{0}^{t}e^{-\delta s}\int\nolimits_{X_{s^{-}}%
}^{X_{s^{-}}+Y_{s^{-}}}\left(  V^{n-1}(X_{s^{-}}+Y_{s^{-}}-\alpha
,0)-\psi(X_{s^{-}},Y_{s^{-}})\right)  dF^{2}(\alpha)ds\\
&  & -\lambda_{2}\int\nolimits_{0}^{t}e^{-\delta s}\int\nolimits_{X_{s^{-}%
}+Y_{s^{-}}}^{\infty}\left(  a_{1}V_{1}^{0}(Y_{s^{-}})-\psi(X_{s^{-}}%
,Y_{s^{-}})\right)  dF^{2}(\alpha)ds\text{. }%
\end{array}
\label{Martingala4}%
\end{equation}
instead of the martingales $M_{t}^{1}$ and $M_{t}^{2\text{ }}$defined in
(\ref{M1t}) and (\ref{M2t}) respectively.\hfill${\small \square}$

In the next proposition we state that $V^{n}$ is the smallest viscosity
solution of the corresponding HJB equation.

\begin{proposition}
\smallskip\label{MenorSuperIterativa} The optimal value function $V^{n}$ is
the smallest viscosity supersolution of (\ref{HJB Iterativa}) satisfying
growth condition A.1.
\end{proposition}

\textit{Proof. }The proof of this proposition is similar to the one of
Proposition \ref{MenorSuper}, but using as martingales with zero expectation
(\ref{Martingala3}) and (\ref{Martingala4}).\hfill{\small $\square$}

\begin{remark}
\label{Verification Iterativa}From the above proposition we deduce the usual
viscosity verification result for the $n-$step: Consider a family of
admissible strategies $\{\overline{L}^{x,y}\in\Pi_{x,y}:\left(  x,y\right)
\in\mathbf{R}_{+}^{2}\}$. If the function $V_{\overline{L}^{x,y}}^{n}(x,y)$ is
a viscosity supersolution of (\ref{HJB Iterativa}) then $V_{\overline{L}%
^{x,y}}^{n}=V^{n}$.
\end{remark}

Finally, we have the convergence result to the optimal value function
(\ref{DefinitionV}).

\begin{proposition}
\label{Convergencia de Vn} $V^{n}\nearrow V$ as $n$ goes to infinity.
\end{proposition}

\textit{Proof. }By Lemma \ref{V increasing_LocallyLip} and Lemma
\ref{V_GrowthCondition}, $V$ is increasing and satisfies property A.1, so
there exists a $T>0$ such that%

\begin{equation}
e^{-\delta t}V(x+p_{1}t,y+p_{2}t)<\frac{\varepsilon}{3}\text{ }
\label{Cota Convergencia1}%
\end{equation}
for $t\geq T$. Let us define $\kappa=V(x+p_{1}T,y+p_{2}T)>0$ and take
$n_{0}>0$ such that
\begin{equation}
P(\tau_{n_{0}}\geq T)\geq1-\frac{\varepsilon}{3\kappa}\text{.}
\label{Cota probabilidad de Tau1}%
\end{equation}
There exists an admissible strategy $\overline{L}\in\Pi_{x,y}$ such that%
\begin{equation}
V(x,y)-V_{\overline{L}}(x,y)\leq\frac{\varepsilon}{3}\text{.}
\label{Desigualdad V y VL}%
\end{equation}
We define the strategy $\overline{L}^{n_{0}}\in\Pi_{x,y}^{n_{0}}$ as
$\overline{L}_{t}^{n_{0}}=\overline{L}_{t}$ for $t\leq\tau_{n_{0}}%
\wedge\overline{\tau}$ and $\overline{L}_{t}^{n_{0}}=\overline{L}%
_{t-\tau_{n_{0}}}^{0}$ for $t\geq\tau_{n_{0}}$ if $\tau_{n_{0}}<\overline
{\tau}$. From (\ref{Cota Convergencia1}), (\ref{Cota probabilidad de Tau1})
and Lemma \ref{V increasing_LocallyLip}, we have%
\[%
\begin{array}
[c]{l}%
V_{\overline{L}}(x,y)-V_{\overline{L}^{n_{0}}}(x,y)\\%
\begin{array}
[c]{ll}%
\leq & E_{x,y}\left(  a_{1}\left(
%TCIMACRO{\tint _{\overline{\tau}\wedge\tau_{n_{0}}}^{\overline{\tau}}}%
%BeginExpansion
{\textstyle\int_{\overline{\tau}\wedge\tau_{n_{0}}}^{\overline{\tau}}}
%EndExpansion
e^{-\delta s}dL_{s}^{1}+e^{-\delta\overline{\tau}}V_{1}^{0}(X_{\overline{\tau
}}^{\overline{L}})\right)  +a_{2}\left(
%TCIMACRO{\tint _{\overline{\tau}\wedge\tau_{n_{0}}}^{\overline{\tau}}}%
%BeginExpansion
{\textstyle\int_{\overline{\tau}\wedge\tau_{n_{0}}}^{\overline{\tau}}}
%EndExpansion
e^{-\delta s}dL_{s}^{2}+e^{-\delta\overline{\tau}}V_{2}^{0}(Y_{\overline{\tau
}}^{\overline{L}})\right)  \right) \\
\leq & E_{x,y}(e^{-\delta\left(  \overline{\tau}\wedge\tau_{n_{0}}\right)
}V(X_{\overline{\tau}\wedge\tau_{n_{0}}}^{\overline{L}},Y_{\overline{\tau
}\wedge\tau_{n_{0}}}^{\overline{L}}))\\
\leq & E_{x,y}(I_{\left\{  \overline{\tau}\wedge\tau_{n_{0}}\geq T\right\}
}e^{-\delta\left(  \overline{\tau}\wedge\tau_{n_{0}}\right)  }V(x+p_{1}\left(
\overline{\tau}\wedge\tau_{n_{0}}\right)  ,y+p_{2}\left(  \overline{\tau
}\wedge\tau_{n_{0}}\right)  )\\
& +E_{x,y}(I_{\left\{  \tau_{n_{0}}<T\right\}  }e^{-\delta\left(
\overline{\tau}\wedge\tau_{n_{0}}\right)  }V(x+p_{1}\left(  \overline{\tau
}\wedge\tau_{n_{0}}\right)  ,y+p_{2}\left(  \overline{\tau}\wedge\tau_{n_{0}%
}\right)  )\\
\leq & E_{x,y}(I_{\left\{  \overline{\tau}\wedge\tau_{n_{0}}\geq T\right\}
}e^{-\delta\left(  \overline{\tau}\wedge\tau_{n_{0}}\right)  }V(x+p_{1}\left(
\overline{\tau}\wedge\tau_{n_{0}}\right)  ,y+p_{2}\left(  \overline{\tau
}\wedge\tau_{n_{0}}\right)  ))+\kappa P(\tau_{n_{0}}<T)\\
\leq & \frac{2\varepsilon}{3}.
\end{array}
\end{array}
\]

Then we obtain from (\ref{Desigualdad V y VL})%

\[
V(x,y)\leq V_{\overline{L}}(x,y)+\frac{\varepsilon}{3}\leq V_{\overline
{L}^{n_{0}}}(x)+\varepsilon\leq V^{n}(x,y)+\varepsilon
\]
for any $n\geq n_{0}$.\hfill{\small $\square$}

\section{Stationary dividend strategies}

\label{sec6}

As in the one-dimensional case (see for instance Azcue and Muler \cite{azmu})
our aim is to find a stationary dividend strategy whose value function is the
optimal value function $V$. A dividend strategy is \textit{stationary} when
the decision on the dividend payment depends on the current surplus only, and
not on the full history of the controlled process; note that a stationary
dividend strategy generates a family of admissible strategies $\left\{
\overline{L}^{x,y}\in\Pi_{x,y}\text{ for any }\left(  x,y\right)
\in\mathbf{R}_{+}^{2}\right\}  $.

If we assume that the optimal value function $V$ is differentiable, the form
in which the optimal value function solves the HJB equation at any
$(x,y)\in\mathbf{R}_{+}^{2}$ suggests how the dividends should be paid when
the current surplus is $(x,y)$. There are only seven possibilities:

\begin{enumerate}
\item[(i)] If the current surplus is in the open set%
\[
\mathcal{C}^{\ast}=\left\{  (x,y)\in\mathbf{R}_{+}^{2}:\mathcal{L}%
(V)(x,y)=0,V_{x}(x,y)>a_{1},V_{y}(x,y)>a_{2}\right\}  ,
\]
no dividends are paid. The set $\mathcal{C}^{\ast}$ is called the
\textit{non-action} set.

\item[(ii)] If the current surplus is in the open set%
\[
\mathcal{B}_{1}^{\ast}=\left\{  (x,y)\in\mathbf{R}_{+}^{2}:\mathcal{L}%
(V)(x,y)<0,V_{x}(x,y)=a_{1},V_{y}(x,y)>a_{2}\right\}  ,
\]
Company One pays a lump sum as dividends. This lump sum should be
$\min\{b>0:(x-b,y)\notin\mathcal{B}_{1}^{\ast}\}.$

\item[(iii)] If the current surplus is in the open set%
\[
\mathcal{B}_{2}^{\ast}=\left\{  (x,y)\in\mathbf{R}_{+}^{2}:\mathcal{L}%
(V)(x,y)<0,V_{x}(x,y)>a_{1},V_{y}(x,y)=a_{2}\right\}  ,
\]
Company Two pays a lump sum as dividends. This lump sum should be
$\min\{b>0:(x,y-b)\notin\mathcal{B}_{2}^{\ast}\}.$

\item[(iv)] If the current surplus is in the set%
\[
\mathcal{B}_{0}^{\ast}=\left\{  (x,y)\in\mathbf{R}_{+}^{2}:\mathcal{L}%
(V)(x,y)<0,V_{x}(x,y)=a_{1},V_{y}(x,y)=a_{2}\right\}  ,
\]
but not in the closure of $\mathcal{B}_{1}^{\ast}\cup\mathcal{B}_{2}^{\ast}$,
both companies pay a lump sum as dividends.

\item[(v)] If the current surplus is in the closed set%
\[
\mathcal{A}_{0}^{\ast}=\left\{  (x,y)\in\mathbf{R}_{+}^{2}:\mathcal{L}%
(V)(x,y)=0,V_{x}(x,y)=a_{1},V_{y}(x,y)=a_{2}\right\}  ,
\]
both companies pay their incoming premiums as dividends.

\item[(vi)] If the current surplus is in the set%
\[
\mathcal{A}_{1}^{\ast}=\left\{  (x,y)\in\mathbf{R}_{+}^{2}:\mathcal{L}%
(V)(x,y)=0,V_{x}(x,y)=a_{1},V_{y}(x,y)>a_{2}\right\}  ,
\]
Company One pays dividends at some special rate so that the surplus remains in
$\mathcal{A}_{1}^{\ast}\cup\mathcal{A}_{0}^{\ast}$.

\item[(vii)] If the current surplus is in the set%
\[
\mathcal{A}_{2}^{\ast}=\left\{  (x,y)\in\mathbf{R}_{+}^{2}:\mathcal{L}%
(V)(x,y)=0,V_{x}(x,y)>a_{1},V_{y}(x,y)=a_{2}\right\}  ,
\]
Company Two pays dividends at some special rate so that the surplus remains in
$\mathcal{A}_{2}^{\ast}\cup\mathcal{A}_{0}^{\ast}$.
\end{enumerate}

Note that if $V$ is a continuously differentiable solution of (\ref{HJB}),
then $\mathcal{A}^{\ast}\mathcal{=\mathcal{A}}_{0}^{\ast}\cup\mathcal{A}%
_{1}^{\ast}\cup\mathcal{A}_{2}^{\ast}$ is closed, $\mathcal{B}^{\ast
}=\mathcal{B}_{1}^{\ast}\cup\mathcal{B}_{2}^{\ast}\cup\mathcal{B}_{0}^{\ast}$
is open, and any segment which connect a point of $\mathcal{B}^{\ast}$ with a
point of the open set $\mathcal{C}^{\ast}$ should contain a point of
$\mathcal{A}^{\ast}$.

\begin{remark}
\label{RemarkSymmetric}Let us consider the simplest case of identical and
independent Cram\'{e}r-Lundberg processes in (\ref{freeSurplus}); that is
$p_{1}=p_{2}=p$; $\lambda_{1}=\lambda_{2}$; and $U_{i}^{1}$, $U_{i}^{2}$ have
the same distribution $F$. We also choose the dividends paid by both
companies to be equally weighted, i.e. $a_{1}=a_{2}=1/2$. Under these
assumptions, the optimal value function will be symmetric, that is
$V(x,y)=V(y,x)$, and so the sets introduced above satisfy the following
properties: the line $y=x$ is an axis of symmetry of the sets $\mathcal{C}%
^{\ast}$, $\mathcal{B}_{0}^{\ast}$ and $\mathcal{A}_{0}^{\ast}$; the sets
$\mathcal{B}_{1}^{\ast}$ and $\mathcal{A}_{1}^{\ast}$ are the reflection with
respect to the line $y=x$ of the sets $\mathcal{B}_{2}^{\ast}$ and
$\mathcal{A}_{2}^{\ast}$ respectively.
\end{remark}

\section{Curve strategies}

\label{sec7}

We introduce a family of stationary dividend strategies, called \textit{curve
strategies}, in which dividends are paid in the seven ways mentioned in the
previous section, having a simple structure: here the boundary between the
action and non-action region is given by a curve. These strategies can be seen
as the natural analogues of the one-dimensional barrier strategies in this
two-dimensional case.

It is reasonable to think that if the optimal strategy is a curve strategy it
should satisfy the following properties: If $(x_{0},y_{0})\in\mathcal{B}%
_{1}^{\ast}\cup\mathcal{A}_{1}^{\ast}$ (that is only Company One pays
dividends), then $(x,y_{0})$ should be in $\mathcal{B}_{1}^{\ast}$ for all
$x>x_{0}$; analogously if $(x_{0},y_{0})\in\mathcal{B}_{2}^{\ast}%
\cup\mathcal{A}_{2}^{\ast}$ (that is only Company Two pays dividends), then
$(x_{0},y)$ should be in $\mathcal{B}_{2}^{\ast}$ for all $y>y_{0}.$ Finally,
the set $\mathcal{C}^{\ast}$ should be bounded because, as in the
one-dimensional case, the surplus of each company under the optimal strategy
should be bounded for $t>0$.

Let us define the curve strategies satisfying the properties mentioned above.
For these strategies, $\mathbf{R}_{+}^{2}$ is partitioned into seven sets
$\mathcal{C}$, $\mathcal{A}_{0}$, $\mathcal{A}_{1}$, $\mathcal{A}_{2}$,
$\mathcal{B}_{0}$, $\mathcal{B}_{1}$ and $\mathcal{B}_{2}$ where
$\mathcal{A}=\mathcal{A}_{0}\cup\mathcal{A}_{1}\cup\mathcal{A}_{2}$ is a curve
which intersects both coordinate axes.

\begin{description}
\item
\begin{itemize}
\item $\mathcal{A}_{0}=\{(\overline{x},\overline{y})\}$ with $\left(
\overline{x},\overline{y}\right)  \in\mathbf{R}_{+}^{2}$. If the current
surplus is $(\overline{x},\overline{y})$, both companies pay their incoming premium as dividends. Let us call $\overline{u}$
the $x$-intercept and $\overline{v}$ the $y$-intercept of the line with slope
$p_{2}/p_{1}$ passing through $(\overline{x},\overline{y});$ let us denote
$\mathcal{O}_{1}^{(\overline{x},\overline{y})}$ and $\mathcal{O}%
_{2}^{(\overline{x},\overline{y})}$ the regions in the first quadrant bounded
above and below by this line, respectively.

\item $\mathcal{B}_{0}=[\overline{x},\infty)\times\lbrack\overline{y}%
,\infty)-\mathcal{A}_{0}$. If the current surplus is $(x,y)\in\mathcal{B}_{0}%
$, Company One and Company Two pay $x-\overline{x}$ and $y-\overline{y}$ as
dividends, respectively.

\item The set $\mathcal{A}_{1}$ is a curve in $\mathcal{O}_{1}^{(\overline
{x},\overline{y})}$ parametrized by
\[
\mathcal{A}_{1}=\left\{  (u+\tfrac{p_{1}}{p_{2}}\xi_{1}(u),\xi_{1}(u))\text{
with }\overline{u}<u\leq M_{\xi_{1}}\right\}  ,
\]
where $\xi_{1}:[\overline{u},M_{\xi_{1}}]\rightarrow\mathbf{R}$ is a
continuously differentiable function with $\xi_{1}(\overline{u})=\overline{y}%
$, $\xi_{1}(M_{\xi_{1}})=0$ and negative derivative. If the current surplus
$(x,y)\in$ $\mathcal{A}_{1}$, Company Two does not pay dividends and Company
One pays dividends at some special rate for which the bivariate surplus
remains in the curve $\mathcal{A}_{1}$. By basic calculus, it can be shown
that this rate is given by
\[
l_{1}(x,y)=-\frac{p_{2}}{\xi_{1}^{\prime}(x-\left(  p_{1}/p_{2}\right)  y)}.
\]

\item The set $\mathcal{B}_{1}$ is the set to the right of $\mathcal{A}_{1}$
in $\mathcal{O}_{1}^{(\overline{x},\overline{y})}$, that is%
\[
\mathcal{B}_{1}=\left\{  (x,y)\in\mathbf{R}_{+}^{2}:y\leq\overline{y}\text{
and }x>\xi_{1}^{-1}(y)+\tfrac{p_{1}}{p_{2}}y\right\}  .
\]
If the current surplus $(x,y)\in$ $\mathcal{B}_{1}$, Company Two does not pay
dividends and Company One pays the lump sum
\[
\min\{b>0:(x-b,y)\in\mathcal{A}_{1}\}=x-(p_{1}/p_{2})y-\xi_{1}%
^{-1}(y).
\]

\item The sets $\mathcal{A}_{2}$ and $\mathcal{B}_{2}$ in $\mathcal{O}%
_{2}^{(\overline{x},\overline{y})}$ are defined analogously to $\mathcal{A}%
_{1}$ and $\mathcal{B}_{1}$ with the roles of Company One and Two
interchanged; that is%
\[
\mathcal{A}_{2}=\left\{  (\xi_{2}(v),v+\tfrac{p_{2}}{p_{1}}\xi_{2}(v))\text{
with }\overline{v}<v\leq M_{\xi_{2}}\right\}  ,
\]
and%
\[
\mathcal{B}_{2}=\left\{  (x,y)\in\mathbf{R}_{+}^{2}:x\leq\overline{x}\text{
and }y>\xi_{2}^{-1}\left(  x\right)  +\tfrac{p_{2}}{p_{1}}x\right\}  ,
\]

where $\xi_{2}:[\overline{v},M_{\xi_{2}}]\rightarrow\mathbf{R}$ is a
continuously differentiable function with $\xi_{2}(\overline{v})=\overline{x}%
$, $\xi_{2}(M_{\xi_{2}})=0$ and negative derivative. If the current surplus
$(x,y)\in$ $\mathcal{A}_{2}$, Company One does not pay dividends and Company
Two pays dividends at some special rate for which the bivariate surplus
remains in the curve $\mathcal{A}_{2}$. Here this rate is
\[
l_{2}(x,y)=-\frac{p_{1}}{\xi_{2}^{\prime}(y-(p_{2}/p_{1})x)}.
\]

\item If the current surplus $(x,y)\in$ $\mathcal{B}_{2}$, Company One does
not pay dividends and Company Two pays the lump sum
\[
\min\{b>0:(x,y-b)\in\mathcal{A}_{2}\}=y-\left(  p_{2}/p_{1}\right)  x-\xi
_{2}^{-1}(x).
\]

\item The no-action region $\mathcal{C}$ is the open set delimited by the
curve $\mathcal{A}$ and the axes. If the current surplus $(x,y)\in$
$\mathcal{C}$, no dividends are paid.
\end{itemize}
\end{description}

The set partition of the curve strategy corresponding to $\left(  \overline
{x},\overline{y}\right)  =(1,2)$ and the functions
\[
\xi_{1}(u)=\tfrac{2(u-4)(u-6)}{35}\text{ for }u\in\left[  -1,4\right]  \text{
and }\xi_{2}(v)=\tfrac{(u-3)(u-6)}{10}\text{ for }v\in\left[  1,3\right]
\]
is illustrated in Figure 7.1.%

\[%
%TCIMACRO{\FRAME{itbpFU}{4.184in}{2.719in}{0in}{\Qcb{Fig. 6.1: Example of Curve
%Strategy}}{}{conjuntospapernew.eps}{\special{ language "Scientific Word";
%type "GRAPHIC";  maintain-aspect-ratio TRUE;  display "USEDEF";
%valid_file "F";  width 4.184in;  height 2.719in;  depth 0in;
%original-width 5.4155in;  original-height 3.506in;  cropleft "0";
%filename 'ConjuntosPaperNew.eps';file-properties "XNPEU";}}}%
%BeginExpansion
{\parbox[b]{4.184in}{\begin{center}
\includegraphics[
height=2.719in,
width=4.184in
]%
{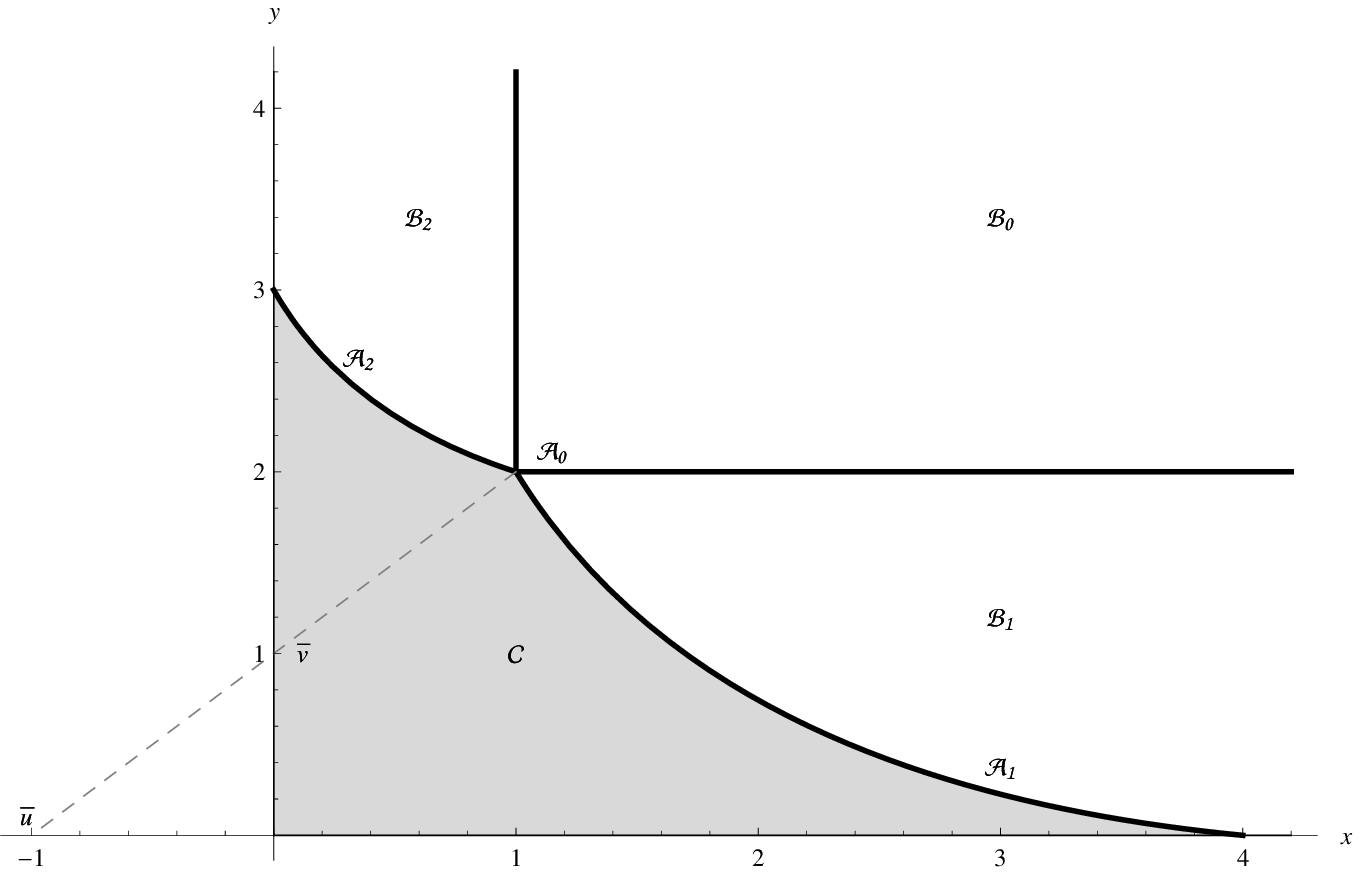}%
\\
Fig. 7.1: Example of Curve Strategy
\end{center}}}%
%EndExpansion
\]

For any $w\in\mathbf{R}$, let us define the set%
\begin{equation}
\Phi^{w}=\left\{  \xi:\left[  w,M_{\xi}\right]  \rightarrow\mathbf{R}%
_{+}\text{, }\xi(M_{\xi})=0,M_{\xi}\geq0\text{, }\xi^{\prime}<0\text{ and }%
\xi^{\prime}\text{ continuous}\right\}  . \label{SetFi}%
\end{equation}

Note that the curve strategies depend only on the point $(\overline
{x},\overline{y})\in\mathbf{R}_{+}^{2}$ and the functions $\xi_{1}\in
\Phi^{\overline{u}}$ and $\xi_{2}\in\Phi^{\overline{v}}$, used in the
parametrization of the curve $\mathcal{A}$. We associate to any $\overline
{\xi}=(\left(  \overline{x},\overline{y}\right)  ,\xi_{1},\xi_{2})$ and any
$(x,y)\in\mathbf{R}_{+}^{2}$ the admissible strategy $\overline{L}%
^{\overline{\xi}}=(L_{t}^{1,\overline{\xi}},L_{t}^{2,\overline{\xi}})\in
\Pi_{x,y}$ Let us define the value function $V^{\overline{\xi}}$ of this curve
strategy as
\begin{equation}
V^{\overline{\xi}}(x,y)=V_{\overline{L}^{\overline{\xi}}}(x,y).
\label{DefinitionV_Psi}%
\end{equation}

We will look for $\overline{\xi}^{\ast}$ such that the associated value
function $V^{\overline{\xi}^{\ast}}$ is the optimal value function defined in
(\ref{DefinitionV}).

\begin{remark}
In the case that $\mathcal{A}$ is the segment $x+y=K$ for some $K>0$ in
$\mathbf{R}_{+}^{2}$, the sum of the dividend rates paid by Company One and
Two is $p_{1}+p_{2}$ for any current surplus in this line. The point
$\mathcal{A}_{0}=(\overline{x},K-\overline{x})$ indicates how this dividend
payment is splitted among the two companies in $\mathcal{A}$: At
$\mathcal{A}_{0}$, Company One pays $p_{1}$ and Company Two pays $p_{2},$ to
the right of this point ($\mathcal{A}_{1}$) Company One pays the total rate
$p_{1}+p_{2}$ and to the left of this point ($\mathcal{A}_{2}$) it is Company Two
which pays $p_{1}+p_{2}$.
\end{remark}

\section{Search for the Optimal Curve Strategy}

\label{sec8}In order to find the optimal value function $V$ of the problem
(\ref{DefinitionV}), we use the iterative approach introduced in Section
\ref{Iterative} and Proposition \ref{Convergencia de Vn}. Our ultimate goal is
to see whether the optimal value function $V$ is the value function of a curve
strategy as defined in the previous section.

We first define an auxiliary function. For any $\overline{\xi}=((\overline
{x},\overline{y}),\xi_{1},\xi_{2}),$ where $(\overline{x},\overline{y}%
)\in\mathbf{R}_{+}^{2}$, $\xi_{1}\in\Phi^{\overline{u}}$ , $\xi_{2}\in
\Phi^{\overline{v}}$ and any continuous function $W_{0}:\mathbf{R}_{+}%
^{2}\rightarrow\lbrack0,+\infty),$ let

\begin{equation}%
\begin{array}
[c]{ll}%
W^{\overline{\xi}\text{ }}(x,y):= & E_{x,y}(\int_{0}^{\tau_{1}}e^{-\delta
s}\left(  a_{1}dL_{s}^{1,\overline{\xi}}+a_{2}dL_{s}^{2,\overline{\xi}%
}\right)  +e^{-\delta\tau_{1}}W_{0}(X_{\tau_{1}}^{\overline{L}^{\xi}}%
,Y_{\tau_{1}}^{\overline{L}^{\xi}})I_{\tau_{1}<\overline{\tau}}\\
& +e^{-\delta\overline{\tau}}\left(  a_{1}V_{1}^{0}(X_{\overline{\tau}%
}^{\overline{L}^{\overline{\xi}}})+a_{2}V_{2}^{0}(Y_{\overline{\tau}%
}^{\overline{L}^{\overline{\xi}}})\right)  I_{\tau_{1}=\overline{\tau}}).
\end{array}
\label{Wpsi}%
\end{equation}

If $W_{0}$ is the value function of a family of admissible strategies
$\overline{L}=\left(  \overline{L}_{x.y}\in\Pi_{x,y}\right)  _{(x,y)\in
\mathbf{R}_{+}^{2}}$, and $\mathcal{C}$, $\mathcal{A}_{0}$, $\mathcal{A}_{1}$,
$\mathcal{A}_{2}$, $\mathcal{B}_{0}$, $\mathcal{B}_{1}$ and $\mathcal{B}_{2}$
are the sets associated to $\overline{\xi}$ (as defined in the previous
section), then $W^{\overline{\xi}\text{ }}$ would be the value function of the
strategy which pays dividends according to the curve strategy $\overline
{L}^{\overline{\xi}}$ up to the first claim and according to $\overline{L}$
afterwards. We call this kind of strategy a \textit{one-step curve strategy}.

Define%

\begin{equation}
H(x,y):=\mathcal{I}(W_{0})(x,y)+U(x,y)\text{.} \label{DefinicionH}%
\end{equation}

In the next proposition, we find an explicit formula for the function
$W^{\overline{\xi}\text{ }}$in terms of $W_{0}$ and $\xi_{1}$ for
$(x,y)\in\mathcal{O}_{1}^{(\overline{x},\overline{y})}$; the formula for the
value function for $(x,y)\in\mathcal{O}_{2}^{(\overline{x},\overline{y})}$
follows in an analogous way and depends only on $\xi_{2}$.

In order to obtain this formula, we use the fact that $W^{\overline{\xi}\text{
}}$satisfies the integro-differential equation $\mathcal{L}^{n}(W^{\overline
{\xi}\text{ }})=0$ in $\mathcal{C\cup A}$ and that $W_{x}^{\overline{\xi
}\text{ }}=a_{1}$ in $\mathcal{A}_{1}\cup\mathcal{A}_{0}\cup\mathcal{B}%
_{1}\cup\mathcal{B}$.

\begin{proposition}
\label{VnFi}Given $\overline{\xi}=((\overline{x},\overline{y}),\xi_{1},\xi
_{2})$ and a continuous function $W_{0},$ we have that%

\[%
\begin{array}
[c]{ll}%
W^{\overline{\xi}}(x,y)= & e^{-(\delta+\lambda)\frac{\xi_{1}(x-\frac{p_{1}%
}{p_{2}}y)-y}{p_{2}}}k(x-\frac{p_{1}}{p_{2}}y)I_{\left\{  \left(
y-\overline{u}\right)  \frac{p_{1}}{p_{2}}\leq x\leq\xi_{1}^{-1}%
(y)+\frac{p_{1}}{p_{2}}y\text{,}y\leq\overline{y}\right\}  }\\
& +(%
%TCIMACRO{\tint _{0}^{\frac{\xi_{1}(x-\frac{p_{1}}{p_{2}}y)-y}{p_{2}}}}%
%BeginExpansion
{\textstyle\int_{0}^{\frac{\xi_{1}(x-\frac{p_{1}}{p_{2}}y)-y}{p_{2}}}}
%EndExpansion
e^{-(\delta+\lambda)w}H(x+p_{1}w,y+p_{2}w)dw)I_{\left\{  \left(
y-\overline{u}\right)  \frac{p_{1}}{p_{2}}\leq x\leq\xi_{1}^{-1}%
(y)+\frac{p_{1}}{p_{2}}y,y\leq\overline{y}\right\}  }\\
& +\left(  a_{1}(x-\xi_{1}^{-1}(y)-\frac{p_{1}}{p_{2}}y)+k(\xi_{1}%
^{-1}(y))\right)  I_{\left\{  x\geq\xi_{1}^{-1}(y)+\frac{p_{1}}{p_{2}}%
y,y\leq\overline{y}\right\}  }\\
& +\left(  a_{1}(x-\overline{x})+a_{2}(y-\overline{y})+k(\overline{x}%
-\frac{p_{1}}{p_{2}}\overline{y})\right)  I_{\left\{  x\geq\overline{x}%
,y\geq\overline{y}\right\}  },
\end{array}
\]

for $(x,y)\in\mathcal{O}_{1}^{(\overline{x},\overline{y})},$ where
$\overline{u}=\overline{x}-\frac{p_{1}}{p_{2}}\overline{y}$, the function $H$
is defined in (\ref{DefinicionH}) and
\end{proposition}

\[%
\begin{array}
[c]{lll}%
k(u) & = & e^{(\delta+\lambda)\frac{\xi_{1}(u)-\xi_{1}(\overline{u})}{p_{2}}%
}\left(  \frac{p}{\delta+\lambda}+\frac{1}{\delta+\lambda}H(\overline{u}%
+\frac{p_{1}}{p_{2}}\xi_{1}(\overline{u}),\xi_{1}(\overline{u}))\right) \\
&  & +e^{(\delta+\lambda)\frac{\xi_{1}(u)}{p_{2}}}a_{1}\int_{\overline{u}}%
^{u}e^{-(\delta+\lambda)\frac{\xi_{1}(w)}{p_{2}}}dw\\
&  & +\frac{e^{(\delta+\lambda)\frac{\xi_{1}(u)}{p_{2}}}}{p_{2}}\int_{\xi
_{1}(u)}^{\xi_{1}(\overline{u})}H(\xi_{1}^{-1}(t)+\frac{p_{1}}{p_{2}%
}t,t)e^{-(\delta+\lambda)\frac{t}{p_{2}}}dt.
\end{array}
\]

\bigskip\textit{Proof. } Let us consider first an initial surplus
$(x,y)\in\mathcal{C}\cap\mathcal{O}_{1}^{(\overline{x},\overline{y})}$. By
definition (\ref{Wpsi}) we have that the controlled surplus process for
$t<\tau_{1}\wedge h$ and $h>0$ small enough is given by%
\[
(X_{t},Y_{t})=(x+p_{1}t,y+p_{2}t).
\]
So we have that%

\[%
\begin{array}
[c]{ll}%
W^{\overline{\xi}}(x,y)= & E_{x,y}(e^{-\delta t}W^{\overline{\xi}}%
(X_{t\wedge\tau_{1}},Y_{t\wedge\tau_{1}})I_{t\wedge\tau_{1}=t}+e^{-\delta
\tau_{1}}W_{0}(X_{\tau_{1}},Y_{\tau_{1}})I_{t\wedge\tau_{1}=\tau_{1}%
<\overline{\tau}}\\
& +e^{-\delta\overline{\tau}}\left(  a_{1}V_{1}^{0}(X_{\overline{\tau}}%
)+a_{2}V_{2}^{0}(Y_{\overline{\tau}})\right)  I_{t\wedge\tau_{1}=\tau
_{1}=\overline{\tau}}).
\end{array}
\]
We can write%

\[%
\begin{array}
[c]{l}%
E_{x,y}\left(  e^{-\delta\,(\tau_{1}\wedge t)}I_{\{\tau_{1}\wedge t=t<\tau
_{1}\}}W^{\overline{\xi}}(X_{\tau_{1}\wedge t}^{\overline{L}},Y_{\tau
_{1}\wedge t}^{\overline{L}})+e^{-\delta\,(\tau_{1}\wedge t)}I_{\{\tau
_{1}\wedge t=\tau_{1}<\overline{\tau}\}}W_{0}(X_{\tau_{1}\wedge t}%
^{\overline{L}},Y_{\tau_{1}\wedge t}^{\overline{L}}))\right) \\%
\begin{array}
[c]{ll}%
= & E_{x,y}\left(  e^{-\delta\,(\tau_{1}\wedge t)}I_{\tau_{1}\wedge
t=t<\tau_{1}}W^{\overline{\xi}}(X_{\tau_{1}\wedge t}^{\overline{L}}%
,Y_{\tau_{1}\wedge t}^{\overline{L}})\right) \\
& +E_{x,y}(I_{\left\{  \tau_{1}=\tau_{1}\wedge t<\overline{\tau}\text{and
}\tau_{1}=\tau_{1}^{1}\right\}  \text{ }}e^{-\delta\tau_{1}^{1}}W_{0}%
(X_{\tau_{1}^{1}}^{\overline{L}},Y_{\tau_{1}^{1}}^{\overline{L}}))\\
& +E_{x,y}(I_{\{\tau_{1}=\tau_{1}\wedge t<\overline{\tau}\text{and }\tau
_{1}=\tau_{1}^{2}\text{ }\}}e^{-\delta\tau_{1}^{2}}W_{0}(X_{\tau_{1}^{2}%
}^{\overline{L}},Y_{\tau_{1}^{2}}^{\overline{L}})),
\end{array}
\end{array}
\]
and so%

\[
\lim_{t\rightarrow0^{+}}\frac{e^{-(\lambda+\delta)t}W^{\overline{\xi}}%
(x+p_{1}t,y+p_{2}t)-W^{\overline{\xi}}(x,y)}{t}=-H(x,y).
\]
Then $g(t)=W^{\overline{\xi}}(x+p_{1}t,y+p_{2}t)$ is continuous and
differentiable as long as $(x+p_{1}t,y+p_{2}t)\in$ $\mathcal{C}$ with%
\begin{equation}
g^{\prime}(0)=(\lambda+\delta)W^{\overline{\xi}}(x,y)-H(x,y).
\label{DerivadaVnC}%
\end{equation}

Let us prove now that the function $W^{\overline{\xi}}$ is continuous in
$\mathcal{A}_{1}$ and has a continuous derivative in the direction of this
curve. In case $(x,y)\in\mathcal{A}_{1}$, we have that for $t<\tau_{1}\wedge
h$ and $h>0$ small enough, the controlled surplus process is
\[
(X_{t},Y_{t})=\Big(x+p_{1}t+{\textstyle\int\nolimits_{0}^{t}}
\frac{p_{2}}{\xi_{1}^{\prime}(X_{s}-\left(  p_{1}/p_{2}\right)  Y_{s})}ds,y+p_{2}t\Big)\in\mathcal{A}_{1}.
\]
By (\ref{Wpsi}), we have that%

\begin{align*}
W^{\overline{\xi}}(x,y)  &  =E_{x,y}\left(  a_{1}\int_{0}^{\tau_{1}\wedge
t}\frac{-p_{2}}{\xi_{1}^{\prime}(X_{s}-\left(  p_{1}/p_{2}\right)  Y_{s}%
)}e^{-\delta s}ds+e^{-\delta t}W^{\overline{\xi}}(X_{t\wedge\tau_{1}%
},Y_{t\wedge\tau_{1}})I_{t\wedge\tau_{1}=t}\right) \\
&  +E_{x,y}\left(  e^{-\delta\tau_{1}}W_{0}(X_{\tau_{1}},Y_{\tau_{1}%
})I_{t\wedge\tau_{1}=\tau_{1}<\overline{\tau}}+e^{-\delta\overline{\tau}%
}\left(  a_{1}V_{1}^{0}(X_{\overline{\tau}})+a_{2}V_{2}^{0}(Y_{\overline{\tau
}})\right)  I_{t\wedge\tau_{1}=\tau_{1}=\overline{\tau}}\right)  .
\end{align*}
Then, with an argument similar to the case of $\mathcal{C}$, we obtain for any
$(x,y)\in\mathcal{A}_{1}$,%

\begin{align*}
&  \lim_{t\rightarrow0}\frac{e^{-(\lambda+\delta)t}W^{\overline{\xi}}\left(
x+p_{1}t+%
%TCIMACRO{\tint \nolimits_{0}^{t}}%
%BeginExpansion
{\textstyle\int\nolimits_{0}^{t}}
%EndExpansion
\frac{p_{2}}{\xi_{1}^{\prime}(X_{s}-\left(  p_{1}/p_{2}\right)  Y_{s}%
)}ds,y+p_{2}t\right)  -W^{\overline{\xi}}(x,y)}{t}\\
&  =-H(x,y)+\frac{a_{1}p_{2}}{\xi_{1}^{\prime}(x-\left(  p_{1}/p_{2}\right)
y)}.
\end{align*}
So
\[
g_{1}(t):=W^{\overline{\xi}}(x+p_{1}t+%
%TCIMACRO{\tint \nolimits_{0}^{t}}%
%BeginExpansion
{\textstyle\int\nolimits_{0}^{t}}
%EndExpansion
\frac{p_{2}}{\xi_{1}^{\prime}(X_{s}-\left(  p_{1}/p_{2}\right)  Y_{s}%
)}ds,y+p_{2}t)
\]
is continuous and differentiable at $t=0$ and satisfies%

\begin{equation}
g_{1}^{\prime}(0)=(\lambda+\delta)W^{\overline{\xi}}(x,y)-H(x,y)+\frac
{a_{1}p_{2}}{\xi_{1}^{\prime}(x-\left(  p_{1}/p_{2}\right)  y)}.
\label{Formula_g1_Prima}%
\end{equation}
Since $(x+p_{1}t+%
%TCIMACRO{\tint \nolimits_{0}^{t}}%
%BeginExpansion
{\textstyle\int\nolimits_{0}^{t}}
%EndExpansion
p_{2}/\xi_{1}^{\prime}(X_{s}-\left(  p_{1}/p_{2}\right)  Y_{s})ds,y+p_{2}%
t)\in\mathcal{A}_{1}$ for $t$ small enough, we have that
\[
(x+p_{1}t+%
%TCIMACRO{\tint \nolimits_{0}^{t}}%
%BeginExpansion
{\textstyle\int\nolimits_{0}^{t}}
%EndExpansion
\frac{p_{2}}{\xi_{1}^{\prime}(X_{s}-\left(  p_{1}/p_{2}\right)  Y_{s}%
)}ds,y+p_{2}t)=(u_{0}(t)+\frac{p_{1}}{p_{2}}\xi_{1}(u_{0}(t)),\xi_{1}%
(u_{0}(t)))
\]
for $u_{0}(t):=x-\left(  p_{1}/p_{2}\right)  y+%
%TCIMACRO{\tint \nolimits_{0}^{t}}%
%BeginExpansion
{\textstyle\int\nolimits_{0}^{t}}
%EndExpansion
p_{2}/\xi_{1}^{\prime}(X_{s}-\left(  p_{1}/p_{2}\right)  Y_{s})ds$; therefore%

\[
W^{\overline{\xi}}(u_{0}(t)+\frac{p_{1}}{p_{2}}\xi_{1}(u_{0}(t)),\xi_{1}%
(u_{0}(t)))=g_{1}(t).
\]
Since $u_{0}^{\prime}(t)=p_{2}/\xi_{1}^{\prime}(X_{t}-\left(  p_{1}%
/p_{2}\right)  Y_{t})$ is continuous and negative, $u_{0}^{-1}$ exists and is
continuously differentiable, so%
\[
k(u):=W^{\overline{\xi}}(u+\frac{p_{1}}{p_{2}}\xi_{1}(u),\xi_{1}%
(u))=g_{1}\circ u_{0}^{-1}(u)
\]
is continuously differentiable.

Defining
\[
W(u,s):=W^{\overline{\xi}}(u+\frac{p_{1}}{p_{2}}\xi_{1}(u)-p_{1}s,\xi
_{1}(u)-p_{2}s),
\]
we obtain for $(u+\left(  p_{1}/p_{2}\right)  \xi_{1}(u)-p_{1}s,\xi
_{1}(u)-p_{2}s)\in\mathcal{C}$ that%

\begin{align*}
&  \widehat{\mathcal{L}}(W)(u,s)\\
&  :=-W_{s}(u,s)-\left(  \delta+\lambda\right)  W(u,s)+H(u+\frac{p_{1}}{p_{2}%
}\xi_{1}(u)-p_{1}s,\xi_{1}(u)-p_{2}s)=0
\end{align*}
The equation $\widehat{\mathcal{L}}(W)(u,s)=0$ is a linear ODE in the variable
$s$, so

\[
W(u,s)e^{(\delta+\lambda)s}-k(u)=%
%TCIMACRO{\tint _{0}^{s}}%
%BeginExpansion
{\textstyle\int_{0}^{s}}
%EndExpansion
e^{(\delta+\lambda)t}H(u+\frac{p_{1}}{p_{2}}\xi_{1}(u)-p_{1}t,\xi_{1}%
(u)-p_{2}t)dt;
\]
therefore%

\begin{align}
&  W^{\overline{\xi}}(u+\frac{p_{1}}{p_{2}}\xi_{1}(u)-p_{1}s,\xi_{1}%
(u)-p_{2}s)\label{Formula1VnPsi}\\
&  =e^{-(\delta+\lambda)s}(k(u)+%
%TCIMACRO{\tint _{0}^{s}}%
%BeginExpansion
{\textstyle\int_{0}^{s}}
%EndExpansion
e^{(\delta+\lambda)t}H(u+\frac{p_{1}}{p_{2}}\xi_{1}(u)-p_{1}t,\xi_{1}%
(u)-p_{2}t)dt),\nonumber
\end{align}
for $\overline{u}\leq u\leq M_{\xi_{1}}$ and $0\leq s\leq\min\left\{  \xi
_{1}(u)/p_{2},u/p_{1}+\xi_{1}(u)/p_{2}\right\}  $. So $W^{\overline{\xi}}$ is
continuously differentiable in the intersection of the set $\mathcal{C}%
\cup\mathcal{A}_{1}$ with $\mathcal{O}_{1}^{(\overline{x},\overline{y})}.$ We
also have, from (\ref{DerivadaVnC}) and (\ref{Formula_g1_Prima}), that for any
$(x,y)\in\mathcal{A}_{1}$,%

\begin{align*}
&  \lim_{t\rightarrow0^{-}}\tfrac{W^{\overline{\xi}}(x+p_{1}t,y+p_{2}%
t)-W^{\overline{\xi}}(x,y)}{t}\\
&  =\lim_{t\rightarrow0^{+}}\tfrac{W^{\overline{\xi}}(x+p_{1}t+\frac{p_{2}%
t}{\xi_{1}^{\prime}(x-\frac{p_{1}}{p_{2}}y)},y+p_{2}t)-W^{\overline{\xi}%
}(x,y)}{t}-a_{1}\tfrac{p_{2}}{\xi_{1}^{\prime}(x-\frac{p_{1}}{p_{2}}y)}.
\end{align*}
Then from%

\begin{align*}
&  p_{1}W_{x^{-}}^{\overline{\xi}.}(x,y)+p_{2}W_{y^{-}}^{\overline{\xi}%
.}(x,y)\\
&  =\left(  p_{1}+\tfrac{p_{2}}{\xi_{1}^{\prime}(x-\frac{p_{1}}{p_{2}}%
y)}\right)  W_{x^{-}}^{\overline{\xi}.}(x,y)+p_{2}W_{y^{-}}^{\overline{\xi}%
.}(x,y)-a_{1}\tfrac{p_{2}}{\xi_{1}^{\prime}(x-\frac{p_{1}}{p_{2}}y)},
\end{align*}
we conclude that $W_{x^{-}}^{\overline{\xi}}(x,y)=a_{1}.$

By (\ref{Formula1VnPsi}), and since $\left(  u+\left(  p_{1}/p_{2}\right)
\xi_{1}(u),\xi_{1}(u)\right)  \in\mathcal{A}_{1}$,%

\[%
\begin{array}
[c]{lll}%
W_{x}^{\overline{\xi}}(u+\frac{p_{1}}{p_{2}}\xi_{1}(u),\xi_{1}(u)) & = &
k^{\prime}(u)+\left(  H(u+\frac{p_{1}}{p_{2}}\xi_{1}(u)),\xi_{1}%
(u))-(\delta+\lambda)k(u)\right)  \tfrac{\xi_{1}^{\prime}(u)}{p_{2}}\\
& = & a_{1},
\end{array}
\]
and then%

\begin{equation}
k(u)=k(\overline{u})e^{(\delta+\lambda)\frac{\xi_{1}(u)-\xi_{1}(\overline{u}%
)}{p_{2}}}+\int_{\overline{u}}^{u}\left(  a_{1}-H(w+\tfrac{p_{1}}{p_{2}}%
\xi_{1}(w)),\xi_{1}(w))\frac{\xi_{1}^{\prime}(w)}{p_{2}}\right)
e^{(\delta+\lambda)\frac{\xi_{1}(u)-\xi_{1}(w)}{p_{2}}}dw. \label{Formula_k_u}%
\end{equation}
At the point $(\overline{u}+\left(  p_{1}/p_{2}\right)  \xi_{1}(\overline
{u}),\xi_{1}(\overline{u}))\in\mathcal{A}_{0}$ the dividend strategy consists
of collecting all the incoming premium as dividends up to the time $\tau_{1}$, so%

\begin{equation}
k(\overline{u})=\frac{p}{\delta+\lambda}+\frac{1}{\delta+\lambda}%
H(\overline{u}+\tfrac{p_{1}}{p_{2}}\xi_{1}(\overline{u}),\xi_{1}(\overline
{u})).\nonumber
\end{equation}
Then we have, from (\ref{Formula_k_u}),%

\[%
\begin{array}
[c]{lll}%
k(u) & = & e^{(\delta+\lambda)\frac{\xi_{1}(u)-\xi_{1}(\overline{u})}{p_{2}}%
}\left(  \frac{p}{\delta+\lambda}+\frac{1}{\delta+\lambda}H(\overline{u}%
+\frac{p_{1}}{p_{2}}\xi_{1}(\overline{u}),\xi_{1}(\overline{u}))\right) \\
&  & +\int_{\overline{u}}^{u}\left(  a_{1}-H(w+\frac{p_{1}}{p_{2}}\xi
_{1}(w)),\xi_{1}(w))\frac{\xi_{1}^{\prime}(w)}{p_{2}}\right)  e^{(\delta
+\lambda)\frac{\xi_{1}(u)-\xi_{1}(w)}{p_{2}}}dw.
\end{array}
\]
We conclude for (\ref{Formula1VnPsi}) that for any $(x,y)$ in the intersection
of the set $\mathcal{C}$ with $\mathcal{O}_{1}^{(\overline{x},\overline{y})}$,%

\[%
\begin{array}
[c]{lll}%
W^{\overline{\xi}}(x,y) & = & e^{-(\delta+\lambda)\frac{\xi_{1}(x-\frac{p_{1}%
}{p_{2}}y)-y}{p_{2}}}k(x-\frac{p_{1}}{p_{2}}y)\\
&  & +%
%TCIMACRO{\tint _{0}^{\frac{\xi_{1}(x-\frac{p_{1}}{p_{2}}y)-y}{p_{2}}}}%
%BeginExpansion
{\textstyle\int_{0}^{\frac{\xi_{1}(x-\frac{p_{1}}{p_{2}}y)-y}{p_{2}}}}
%EndExpansion
e^{-(\delta+\lambda)w}H(x+p_{1}w,y+p_{2}w)dw,
\end{array}
\]
which yields the result.\hfill  {\small $\square$}

\begin{remark}
\label{Formula_VPsi_O2}The formula of $W^{\overline{\xi}}$ in $\mathcal{O}%
_{2}^{(\overline{x},\overline{y})}$ can be obtained using the formula given in
Proposition \ref{VnFi} by interchanging the role of Company One and Company
Two using that $W_{y}^{\overline{\xi}\text{ }}=a_{2}$ in $\mathcal{A}_{2}%
\cup\mathcal{A}_{0}\cup\mathcal{B}_{2}\cup\mathcal{B}_{0}$. More precisely, if
$(x,y)\in\mathcal{O}_{2}^{(\overline{x},\overline{y})},$%

\[%
\begin{array}
[c]{ll}%
W^{\overline{\xi}}(x,y)= & e^{-(\delta+\lambda)\frac{\xi_{2}(y-\frac{p_{2}%
}{p_{1}}x)-x}{p_{1}}}\widetilde{k}(y-\frac{p_{2}}{p_{1}}x)I_{\left\{
x\leq\overline{x},y\leq\xi_{2}^{-1}(x)+\frac{p_{2}}{p_{1}}x\right\}  }\\
& +(%
%TCIMACRO{\tint _{0}^{\frac{\xi_{2}(y-\frac{p_{2}}{p_{1}}x)-x}{p_{1}}}}%
%BeginExpansion
{\textstyle\int_{0}^{\frac{\xi_{2}(y-\frac{p_{2}}{p_{1}}x)-x}{p_{1}}}}
%EndExpansion
e^{-(\delta+\lambda)w}H(x+p_{1}w,y+p_{2}w)dw)I_{\left\{  x\leq\overline
{x},y\leq\xi_{2}^{-1}(x)+\frac{p_{2}}{p_{1}}x\right\}  }\\
& +\left(  a_{2}(y-\xi_{2}^{-1}(x)-\frac{p_{2}}{p_{1}}x)+\widetilde{k}(\xi
_{2}^{-1}(x))\right)  I_{\left\{  x\leq\overline{x},y\geq\xi_{2}^{-1}%
(x)+\frac{p_{2}}{p_{1}}x\right\}  }\\
& +\left(  a_{1}(x-\overline{x})+a_{2}(y-\overline{y})+\widetilde{k}%
(\overline{y}-\frac{p_{2}}{p_{1}}\overline{x})\right)  I_{\left\{
x\geq\overline{x},y\geq\overline{y}\right\}  },
\end{array}
\]
where%
\[%
\begin{array}
[c]{lll}%
\widetilde{k}(v) & = & e^{(\delta+\lambda)\frac{\xi_{2}(v)-\xi_{2}%
(\overline{v})}{p_{1}}}\left(  \frac{p}{\delta+\lambda}+\frac{1}%
{\delta+\lambda}H(\xi_{2}(\overline{v}),\overline{v}+\frac{p_{2}}{p_{1}}%
\xi_{2}(\overline{v}))\right) \\
&  & +e^{(\delta+\lambda)\frac{\xi_{2}(v)}{p_{1}}}a_{2}\int_{\overline{v}}%
^{v}e^{-(\delta+\lambda)\frac{\xi_{2}(w)}{p_{1}}}dw\\
&  & +\frac{e^{(\delta+\lambda)\frac{\xi_{2}(v)}{p_{1}}}}{p_{1}}\int_{\xi
_{2}(v)}^{\xi_{2}(\overline{v})}H(t,\xi_{2}^{-1}(t)+\frac{p_{2}}{p_{1}%
}t)e^{-(\delta+\lambda)\frac{t}{p_{1}}}dt.
\end{array}
\]

\end{remark}

From the formulas obtained in Proposition \ref{VnFi} and Remark
\ref{Formula_VPsi_O2} we obtain the following regularity result.

\begin{proposition}
If the function $H$ defined in (\ref{DefinicionH}) is continuously
differentiable, then $W^{\overline{\xi}}$ is continuously differentiable in
$\mathbf{R}_{+}^{2}$.
\end{proposition}

\bigskip\textit{Proof. }Since $\xi_{1}$ and $\xi_{2}$ are continuously differentiable, it is
clear that $W^{\overline{\xi}}$ is continuously differentiable except possibly
at the points of either the boundary of $\mathcal{B}_{0}$ or the segment
\[
\mathcal{S}=\left\{  (x,\tfrac{p_{2}}{p_{1}}(x-\overline{x})+\overline{y}%
)\in\mathbf{R}_{+}^{2}\text{ with }x\leq\overline{x}\right\}  .
\]
After some easy calculations and using that $W^{\overline{\xi}}$ satisfies%
\[
p-(\delta+\lambda)W^{\overline{\xi}}(\overline{x},\overline{y})+H(\overline
{x},\overline{y})=0,\text{ }%
\]
it can be seen that $W^{\overline{\xi}}$ is continuously differentiable in
$\mathcal{S}$ with%

\[
W_{x}^{\overline{\xi}}(x,\tfrac{p_{2}}{p_{1}}(x-\overline{x})+\overline{y})=%
%TCIMACRO{\tint _{0}^{\frac{\overline{x}-x}{p_{1}}}}%
%BeginExpansion
{\textstyle\int_{0}^{\frac{\overline{x}-x}{p_{1}}}}
%EndExpansion
e^{-(\delta+\lambda)w}H_{x}(x+p_{1}w,\tfrac{p_{2}}{p_{1}}(x-\overline
{x})+\overline{y}+p_{2}w)dw+e^{-(\delta+\lambda)\frac{\overline{x}-x}{p_{1}}%
}a_{1}%
\]
and%

\[
W_{y}^{\overline{\xi}}(x,\tfrac{p_{2}}{p_{1}}(x-\overline{u}))=%
%TCIMACRO{\tint _{0}^{\frac{\overline{x}-x}{p_{1}}}}%
%BeginExpansion
{\textstyle\int_{0}^{\frac{\overline{x}-x}{p_{1}}}}
%EndExpansion
e^{-(\delta+\lambda)w}H_{y}(x+p_{1}w,\tfrac{p_{2}}{p_{1}}(x-\overline
{x})+\overline{y}+p_{2}w)dw+e^{-(\delta+\lambda)\frac{\overline{x}-x}{p_{1}}%
}a_{2}.
\]
Finally, the differentiability at the boundary of $\mathcal{B}_{0}$ follows
from the differentiability of $W^{\overline{\xi}}$ at $(\overline{x}%
,\overline{y})$ of $\mathcal{S}$.\hfill {\small $\square$}

Let us define the set of functions
\[
\mathcal{M}=\{w:\mathbf{R}_{+}^{2}\rightarrow\lbrack0,+\infty)\text{
continuous with }w(x,y)-a_{1}x-a_{2}y\text{ bounded}\}.
\]

\begin{proposition}
\label{PuntoFijo_Vpsi} The value function $V^{\overline{\xi}}$ of the curve
strategy corresponding to $\overline{\xi}=((\overline{x},\overline{y}),\xi
_{1},\xi_{2})$ as defined in (\ref{DefinitionV_Psi}), satisfies the formulas
given in Proposition \ref{VnFi} and Remark \ref{Formula_VPsi_O2} replacing
both $W_{0}$ and $W^{\overline{\xi}}$ by $V^{\overline{\xi}}$. Moreover, $V^{\overline{\xi}}$ is
the unique function in $\mathcal{M}$ which satisfies this property.
\end{proposition}

\bigskip\textit{Proof. } $\mathcal{M}$ is a complete metric space with the distance
$d(w_{1},w_{2})=\sup_{\mathbf{R}_{+}^{2}}\left\vert w_{1}-w_{2}\right\vert .$
The operator $\mathcal{T}:\mathcal{M}\rightarrow\mathcal{M}$ defined as%

\[%
\begin{array}
[c]{ll}%
\mathcal{T}(w)(x,y):= & E_{x,y}(\int_{0}^{\tau_{1}}e^{-\delta s}\left(
a_{1}dL_{s}^{1,\overline{\xi}}+a_{2}dL_{s}^{2,\overline{\xi}}\right)
+e^{-\delta\tau_{1}}w(X_{\tau_{1}}^{\overline{L}^{\xi}},Y_{\tau_{1}%
}^{\overline{L}^{\xi}})I_{\tau_{1}<\overline{\tau}}\\
& +e^{-\delta\overline{\tau}}\left(  a_{1}V_{1}^{0}(X_{\overline{\tau}%
}^{\overline{L}^{\overline{\xi}}})+a_{2}V_{2}^{0}(Y_{\overline{\tau}%
}^{\overline{L}^{\overline{\xi}}})\right)  I_{\tau_{1}=\overline{\tau}})
\end{array}
\]
is a contraction with contraction factor $\lambda/(\delta+\lambda)<1.$ Then,
there exists a unique fixed point and by definition (\ref{DefinitionV_Psi}) ,
$\mathcal{T}(V^{\overline{\xi}})=V^{\overline{\xi}}$. Taking in Proposition
\ref{VnFi} and in Remark \ref{Formula_VPsi_O2} the function $W_{0}$ as
$V^{\overline{\xi}}$ we obtain from (\ref{Wpsi}) that $V^{\overline{\xi}%
}=W^{\overline{\xi}}$ and so we get the result.\hfill {\small $\square$}\\

This last proposition gives a constructive way to obtain $V^{\overline{\xi}}$.
Starting with $w_{0}(x,y)=a_{1}x+a_{2}y\in\mathcal{M}$, we define iteratively
$w_{n+1}=\mathcal{T}(w_{n}).$ Hence, $V^{\overline{\xi}}=\lim_{n\rightarrow
\infty}w_{n}$. Note that at each step $w_{n+1}$ can be obtained from the
formulas given in Proposition \ref{VnFi} and Remark \ref{Formula_VPsi_O2}
replacing $W_{0}$ by $w_{n}$.

Consider now the function $V^{n,\overline{\xi}\text{ }}$ defined in
(\ref{Wpsi}) taking $W_{0}$ as the optimal value function $V^{n-1}$
corresponding to step $n-1$ in (\ref{Definicion Vn}). We try to find
$\overline{\xi}_{n}^{\ast}$, which maximizes $V^{n,\overline{\xi}\text{ }}%
$among all the possible $\overline{\xi}=((\overline{x},\overline{y}),\xi
_{1},\xi_{2})$. If the function $V^{n,\overline{\xi}_{n}^{\ast}\text{ }}$ is a
viscosity supersolution of (\ref{HJB Iterativa}), then by Remark
\ref{Verification Iterativa}, we would have that $V^{n,\overline{\xi}%
_{n}^{\ast}\text{ }}=V^{n}$. In the case that one-step curve strategies
corresponding to $\overline{\xi}_{n}^{\ast}$ exist for all $n\geq1$, by
Proposition \ref{Convergencia de Vn}, $V^{n,\overline{\xi}_{n}^{\ast}\text{ }%
}\nearrow V$.

Let us call, as in (\ref{DefinicionH}),%

\begin{equation}
H_{n-1}(x,y):=\mathcal{I}(V^{n-1})(x,y)+U(x,y)\text{ }. \label{Def_H_n-1}%
\end{equation}
In order to find the optimal one-step curve strategy corresponding to
$\overline{\xi}_{n}^{\ast}=((\overline{x}_{n}^{\ast},\overline{y}_{n}^{\ast
}),\xi_{1,n}^{\ast},\xi_{2,n}^{\ast})$, we look first for the optimal vertex
$(\overline{x}_{n}^{\ast},\overline{y}_{n}^{\ast})$. By the formula given in
Proposition \ref{VnFi},%

\[
V^{n,\overline{\xi}}(\overline{x},\overline{y})=\frac{p}{\delta+\lambda}%
+\frac{H_{n-1}(\overline{x},\overline{y})}{\delta+\lambda}%
\]
and
\begin{align*}
&  V^{n,\overline{\xi}}(x,y)\\
&  =V^{n,\overline{\xi}}(\overline{x},\overline{y})+a_{1}(x-\overline
{x})+a_{2}(y-\overline{y}).
\end{align*}
for $x$ and $y$ large enough. So%

\[
(\overline{x}_{n}^{\ast},\overline{y}_{n}^{\ast})=\arg\max_{(\overline
{x},\overline{y})\in\mathbf{R}_{+}^{2}}\frac{H_{n-1}(\overline{x},\overline
{y})}{\delta+\lambda}-a_{1}\overline{x}-a_{2}\overline{y}.
\]
If this maximum is attained at a critical point (assuming that $H_{n-1}$ is
differentiable), we have that $(\overline{x}_{n}^{\ast},\overline{y}_{n}%
^{\ast})$ is a solution of%

\[
\left\{
\begin{array}
[c]{l}%
\partial_{x}H_{n-1}(\overline{x},\overline{y})=a_{1}\left(  \delta
+\lambda\right) \\
\partial_{y}H_{n-1}(\overline{x},\overline{y})=a_{2}\left(  \delta
+\lambda\right)  .
\end{array}
\right.
\]
Let us call $\overline{u}_{n}^{\ast}=\overline{x}_{n}^{\ast}-\left(
p_{1}/p_{2}\right)  \overline{y}_{n}^{\ast}$ and $\overline{v}_{n}^{\ast
}=\overline{y}_{n}^{\ast}-\left(  p_{2}/p_{1}\right)  \overline{x}_{n}^{\ast
}.$ Next, we use Calculus of Variations in order to find two curves $\xi
_{1,n}^{\ast}$ and $\xi_{2,n}^{\ast}$ which maximize $V^{n,\overline{\xi
}\text{ }}(x,y)$, among all $\overline{\xi}=((\overline{x}_{n}^{\ast
},\overline{y}_{n}^{\ast}),\xi_{1},\xi_{2})$ for $\xi_{1}\in\Phi^{\overline
{u}_{n}^{\ast}}$, $\xi_{2}\in\Phi^{\overline{v}_{n}^{\ast}}$ and $\left(
x,y\right)  $ large enough. The two curves can be obtained separately and independently.

\begin{proposition}
\label{E-L-Equation} Assume that $H_{n-1}$ is differentiable and that there
exists $\overline{\xi}_{n}^{\ast}=((\overline{x}_{n}^{\ast},\overline{y}%
_{n}^{\ast}),\xi_{1,n}^{\ast},\xi_{2,n}^{\ast})$ where $\xi_{1,n}^{\ast}%
\in\Phi^{\overline{u}_{n}^{\ast}}$ and $\xi_{2,n}^{\ast}\in\Phi^{\overline
{v}_{n}^{\ast}}$ such that $V^{n}=V^{n,\overline{\xi}_{n}^{\ast}\text{ }}$.
Then $\xi_{1,n}^{\ast}$ satisfies%
\[
\partial_{x}H_{n-1}(u+\tfrac{p_{1}}{p_{2}}\xi_{1,n}^{\ast}(u),\xi_{1,n}^{\ast
}(u))=a_{1}(\delta+\lambda)
\]
for $\overline{u}_{n}^{\ast}\leq u\leq M_{\xi_{1,n}^{\ast}},$ and $\xi
_{2,n}^{\ast}$ satisfies%
\[
\partial_{y}H_{n-1}(\xi_{2,n}^{\ast}(v),v+\tfrac{p_{2}}{p_{1}}\xi_{2,n}^{\ast
}(v))=a_{2}(\delta+\lambda)
\]
for $\overline{v}_{n}\leq v\leq M_{\xi_{2,n}^{\ast}}.$
\end{proposition}

\bigskip\textit{Proof. }
We will prove this result for $\xi_{1,n}^{\ast}$, the proof for $\xi
_{2,n}^{\ast}$ is analogous.

Given any $\xi_{1}\in\Phi^{\overline{u}_{n}^{\ast}},$ we have that
\[
V^{n,\overline{\xi}}(M_{\xi_{1}},0)+a_{1}(x-M_{\xi_{1}}),
\]
for $\overline{\xi}=((\overline{x}_{n}^{\ast},\overline{y}_{n}^{\ast}),\xi
_{1},\xi_{2})$ and $x\geq M_{\xi_{1}}$. Then, if there exists $\xi_{1,n}%
^{\ast}\in\Phi^{\overline{u}_{n}^{\ast}}$ such that $V^{n}=V^{n,\overline{\xi
}_{n}^{\ast}}$,%

\[
V^{n,\overline{\xi}_{n}^{\ast}}(M_{\xi_{1,n}^{\ast}},0)-a_{1}M_{\xi
_{1,n}^{\ast}}=\max_{\xi_{1}\in\Phi^{\overline{u}_{n}^{\ast}}}\left(
V^{n,\overline{\xi}}(M_{\xi_{1}},0)-a_{1}M_{\xi_{1}}\right)  .
\]
Consider non-negative test functions $\varsigma$ with $\varsigma(\overline
{u}_{n}^{\ast})=0$ and $\varsigma(M_{\xi_{1,n}^{\ast}})$ $=0.$ We have that
$\xi_{1,n}^{\ast}+\varepsilon\varsigma\in\Phi^{\overline{u}_{n}^{\ast}}$ for
$\varepsilon$ small enough$.$ Let us write,%
\[
\xi_{\varepsilon}(u)=\xi_{1,n}^{\ast}(u)+\varepsilon\varsigma(u).
\]
We have that $M_{\xi_{\varepsilon}}=M_{\xi_{1,n}^{\ast}}$ and then
\[
V^{n,\overline{\xi}_{n}^{\ast}}(M_{\xi_{1,n}^{\ast}},0)-a_{1}M_{\xi
_{1,n}^{\ast}}=\max_{\varsigma}\left(  V^{n,\overline{\xi}_{\varepsilon}%
}(M_{\xi_{1,n}^{\ast}},0)-a_{1}M_{\xi_{1,n}^{\ast}}\right)  ,
\]
where $\overline{\xi}_{\varepsilon}=((\overline{x}_{n}^{\ast},\overline{y}%
_{n}^{\ast}),\xi_{\varepsilon},\xi_{2}).$ Denote%

\[
\vartheta(\varepsilon):=V^{n,\overline{\xi}_{\varepsilon}}(M_{\xi_{1,n}^{\ast
}},0)-a_{1}M_{\xi_{1,n}^{\ast}}.
\]
We have that $\xi_{\varepsilon}(M_{\xi_{1,n}^{\ast}})=\xi_{1,n}^{\ast}%
(M_{\xi_{1,n}^{\ast}})=0$ and $\xi_{\varepsilon}(\overline{u}_{n}^{\ast}%
)=\xi_{1,n}^{\ast}(\overline{u}_{n}^{\ast})$ so we can write by Proposition
\ref{VnFi},

\begin{align}
\vartheta(\varepsilon)  &  =e^{-\frac{(\delta+\lambda)\xi_{1,n}^{\ast
}(\overline{u}_{n}^{\ast})}{p_{2}}}\left(  \frac{p}{\delta+\lambda}+\frac
{1}{\delta+\lambda}H_{n-1}(\overline{u}_{n}^{\ast}+\tfrac{p_{1}}{p_{2}}%
\xi_{1,n}^{\ast}(\overline{u}_{n}^{\ast}),\xi_{1,n}^{\ast}(\overline{u}%
_{n}^{\ast}))\right)  -a_{1}\overline{u}_{n}^{\ast}\nonumber\\
&  +\int_{\overline{u}_{n}^{\ast}}^{M_{\xi_{1,n}^{\ast}}}a_{1}\left(
e^{-\frac{(\delta+\lambda)\xi_{\varepsilon}(w)}{p_{2}}}-1\right)  dw+\frac
{1}{p_{2}}\int_{0}^{\xi_{1,n}^{\ast}(\overline{u}_{n}^{\ast})}H_{n-1}%
(\xi_{\varepsilon}^{-1}(t)+\tfrac{p_{1}}{p_{2}}t,t)e^{-\frac{(\delta
+\lambda)t}{p_{2}}}dt.\nonumber
\end{align}
Clearly,%

\begin{align*}
0  &  =\left.  \frac{\vartheta(\varepsilon)}{\partial\varepsilon}\right\vert
_{\varepsilon=0}\\
&  =\frac{1}{p_{2}}\int_{\overline{u}_{n}^{\ast}}^{M_{\xi_{1,n}^{\ast}}%
}\left(  \partial_{x}H_{n-1}(w+\tfrac{p_{1}}{p_{2}}\xi_{1,n}^{\ast}%
(w),\xi_{1,n}^{\ast}(w))-a_{1}(\delta+\lambda)\right)  e^{-(\delta
+\lambda)\frac{\xi_{1,n}^{\ast}(w)}{p_{2}}}\varsigma(w)dw.
\end{align*}
So we obtain%

\[
\left(  \partial_{x}H_{n-1}(u+\tfrac{p_{1}}{p_{2}}\xi_{1,n}^{\ast}%
(u),\xi_{1,n}^{\ast}(u))-a_{1}(\delta+\lambda)\right)  e^{-(\delta
+\lambda)\frac{\xi_{1,n}^{\ast}(u)}{p_{2}}}=0
\]
for all $\overline{u}_{n}^{\ast}\leq u\leq M_{\xi_{1,n}^{\ast}}.$\hfill ${\small \square}$

This last proposition gives us a constructive way to find the candidate for
$V^{n}$ in the case that it comes from a one-step curve strategy$.$ We find
numerically, if it exists, the solution $z_{1}(u)$ of the equation%

\[
\partial_{x}H_{n-1}(u+\tfrac{p_{1}}{p_{2}}z_{1}(u),z_{1}(u))=a_{1}%
(\delta+\lambda)
\]
for $\overline{u}_{n}^{\ast}\leq u\leq\min\left\{  u:z_{1}(u)=0\right\}  $ and
the solution $z_{2}(v)$ of the equation%
\[
\partial_{y}H_{n-1}(z_{2}(v),v+\tfrac{p_{2}}{p_{1}}z_{2}(v))=a_{2}%
(\delta+\lambda)
\]
for $\overline{v}_{n}^{\ast}\leq v\leq\min\left\{  v:z_{2}(v)=0\right\}  $ .
If $z_{1}(u)$ is in $\Phi^{\overline{u}_{n}^{\ast}}$ and $z_{2}(v)$ is in
$\Phi^{\overline{v}_{n}}$, we define $\xi_{1,n}^{\ast}(u)=z_{1}(u)$ and
$\xi_{2,n}^{\ast}(v)=z_{2}(v)$ and we obtain the value function
$V^{n,\overline{\xi}_{n}^{\ast}}$ by the formula given in Proposition
\ref{VnFi}; this is our candidate for $V^{n}$. Afterwards, we check whether
$V^{n,\overline{\xi}_{n}^{\ast}}$ is a viscosity supersolution of
(\ref{HJB Iterativa}); if this is the case, then $V^{n}=V^{n,\overline{\xi
}_{n}^{\ast}}$.\\

In the next Proposition we state some conditions under which the optimal
strategy of (\ref{DefinitionV}) is a curve strategy. This result, together
with Propositions \ref{VnFi} and \ref{E-L-Equation} gives a way to find the
optimal curve (if it exists). Let us first define a criterion of convergence
for a sequence $\left(  \overline{\xi}_{n}\right)  _{n\geq1}$ that will be
used in the next proposition.

\begin{definition}
\label{Convergence Psik}We say that $\overline{\xi}_{n}=((\overline{x}%
_{n},\overline{y}_{n}),\xi_{1,n},\xi_{2,n})$ converges to $\overline{\xi
}=((\overline{x},\overline{y}),\xi_{1},\xi_{2})$ if
\end{definition}

\[%
\begin{array}
[c]{l}%
\lim\limits_{n\rightarrow\infty}(\overline{x}_{n},\overline{y}_{n}%
)=(\overline{x},\overline{y})\text{, }\lim\limits_{n\rightarrow\infty}%
M_{\xi_{i,n}}\rightarrow M_{\xi_{i}}\text{ for }i=1,2,\\
\lim\limits_{n\rightarrow\infty}\max\limits_{{\small [}\overline{u}%
_{n}{\small ,M}_{\xi_{1,n}}{\small ]\cap\lbrack}\overline{u}{\small ,M}%
_{\xi_{1}}{\small ]}}\left\vert \xi_{1,n}(u)-\xi_{1}(u)\right\vert =0\text{,
}\lim\limits_{n\rightarrow\infty}\max\limits_{[\overline{v}_{n}{\small ,M}%
_{\xi_{2,n}}]\cap\lbrack\overline{v}{\small ,M}_{\xi_{2}}]}\left\vert
\xi_{2,n}(v)-\xi_{2}(v)\right\vert =0,\\
\lim\limits_{n\rightarrow\infty}\max\limits_{[{\small 0,\xi}_{1,n}%
{\small (}\overline{u}_{n}{\small )}]\cap\lbrack{\small 0,\xi}_{1}%
{\small (}\overline{u}{\small )}]}\left\vert \xi_{1,n}^{-1}(w)-\xi_{1}%
^{-1}(w)\right\vert =0\text{ and }\lim\limits_{n\rightarrow\infty}%
\max\limits_{[{\small 0,\xi}_{2,n}{\small (}\overline{v}_{n}{\small )}%
]\cap\lbrack{\small 0,\xi}_{2}{\small (}\overline{v}{\small )}]}\left\vert
\xi_{2,n}^{-1}(w)-\xi_{2}^{-1}(w)\right\vert =0.
\end{array}
\]

\begin{proposition}
\label{ConvergenceCurvan_a_Optima}Assume that there exists a $\overline{\xi}%
_{n}^{\ast}$ such that $V^{n}=V^{n,\overline{\xi}_{n}^{\ast}}$ for all
$n\geq1.$ If $\overline{\xi}_{n}^{\ast}$ converges to some $\overline{\xi
}^{\ast}$ in the sense of Definition \ref{Convergence Psik}, then the optimal
value function $V$ is the value function of the curve strategy $V^{\overline
{\xi}^{\ast}}$ as defined in (\ref{DefinitionV_Psi}).
\end{proposition}

\bigskip\textit{Proof. }
From Proposition \ref{Convergencia de Vn}, we have that $\lim_{n\rightarrow
\infty}V^{n,\overline{\xi}_{n}^{\ast}}=V$. So replacing $W_{0}$ by
$V^{n,\overline{\xi}_{n}^{\ast}}$ and $\overline{\xi}$ by $\overline{\xi}%
_{n}^{\ast}$ in the formulas given in Proposition \ref{VnFi} and Remark
\ref{Formula_VPsi_O2} and letting $n$ go to infinity, we obtain that $V$
satisfies the formulas given in Proposition \ref{VnFi} and Remark
\ref{Formula_VPsi_O2} replacing $W_{0}$ by $V$ and $\overline{\xi}$ by
$\overline{\xi}^{\ast}$. Therefore, by Proposition \ref{PuntoFijo_Vpsi} and
Lemma \ref{V_GrowthCondition}, the functions $V$ and $V^{\overline{\xi}^{\ast
}}$ coincide.\hfill {\small $\square$}

\section{Numerical Example}

\label{sec9}

We present a numerical example in the symmetric and equally weighted case with
an exponential claim size distribution. By Remark \ref{RemarkSymmetric}, we
restrict the search of the optimal curve strategy to $\overline{\xi}=(\left(
\overline{x},\overline{x}\right)  ,\xi,\xi),$ $\xi\in\Phi^{0}$ . Using the
formulas given in Propositions \ref{VnFi} and \ref{E-L-Equation}, we obtain
the functions $\xi_{n}^{\ast}\in\Phi^{0}$ and we check numerically that for
$\overline{\xi}_{n}^{\ast}=(\left(  \overline{x}_{n}^{\ast},\overline{x}%
_{n}^{\ast}\right)  ,\xi_{n}^{\ast},\xi_{n}^{\ast})$ the associated value
function $V^{n,\overline{\xi}_{n}^{\ast}}$ is a viscosity solution of
(\ref{HJB Iterativa}). We also obtain numerically the convergence of
$\overline{\xi}_{n}^{\ast}$ to $\overline{\xi}$ according to Definition
\ref{Convergence Psik}. Then, using Proposition
\ref{ConvergenceCurvan_a_Optima}, one can conclude that the optimal strategy
is a curve strategy with curve $\xi^{\ast}$.

The numerical procedure was done with the \textit{Mathematica} software and
the calculation is quite time-consuming. The concrete chosen parameters are:
exponential claim size distribution with parameter $3$, Poisson intensity
$\lambda_{1}=\lambda_{2}=20/9$, premium rate $p_{1}=p_{2}=1$, and a discount
factor $\delta=0.1$. In this numerical procedure we used step-size $\Delta
x=\Delta y=0.002$ and iterated $60$ times. The resulting optimal curve
strategy is given in Figure 9.1, and $V(x,y)-(x+y)/2$ (the improvement
of the optimal dividend strategy over paying out the initial capital
immediately) is depicted as the upper curve in Figure 9.2.\newline

We also compare for this numerical example the optimal value function $V(x,y)$
with the (comparably weighted) sum of the stand-alone value functions without
collaboration:%
\[
V_{S}(x,y)=\frac{V^{0}(x)+V^{0}(y)}{2}%
\]
and with $V_{M}(x+y)/2,$ where $V_{M}$ is the optimal value function for the
merger of the two companies. Figure 9.2 depicts the graphics of all
three value functions $V(x,y)$, $V_{S}(x,y)$ and $V_{M}(x+y)/2,$ each of them
reduced by $(x+y)/2$. The optimal merger strategy is barrier with barrier
$b=2.77.$ By Remark \ref{Comparacion Merger}, $V_{M}(x+y)/2<V(x,y)$ for all
$(x,y)\in\mathbf{R}_{+}^{2}$. One sees that whereas for the comparison between
the stand-alone case and the merger the initial surplus levels matter (with
the merger case being the lowest of the three value functions in $(0,0)$), the
collaboration case outperforms not only the merger case but also the
stand-alone one for all combinations of initial surplus levels (i.e. if one
measures the overall dividend payments that can be achieved with either
behavior, for this numerical example collaboration is always preferable).
Hence we have here an instance where collaboration is beneficial not only for
safety aspects, but also with respect to collective profitability.%

\[%
%TCIMACRO{\FRAME{itbpFU}{2.4829in}{2.5339in}{0in}{\Qcb{Figure 9.1: Optimal
%Curve Strategy}}{}{optimopapernew.eps}{\special{ language "Scientific Word";
%type "GRAPHIC";  maintain-aspect-ratio TRUE;  display "USEDEF";
%valid_file "F";  width 2.4829in;  height 2.5339in;  depth 0in;
%original-width 2.7596in;  original-height 2.8167in;  cropleft "0";
%filename 'OptimoPaperNew.eps';file-properties "XNPEU";}}}%
%BeginExpansion
{\parbox[b]{2.4829in}{\begin{center}
\includegraphics[
height=2.5339in,
width=2.4829in
]%
{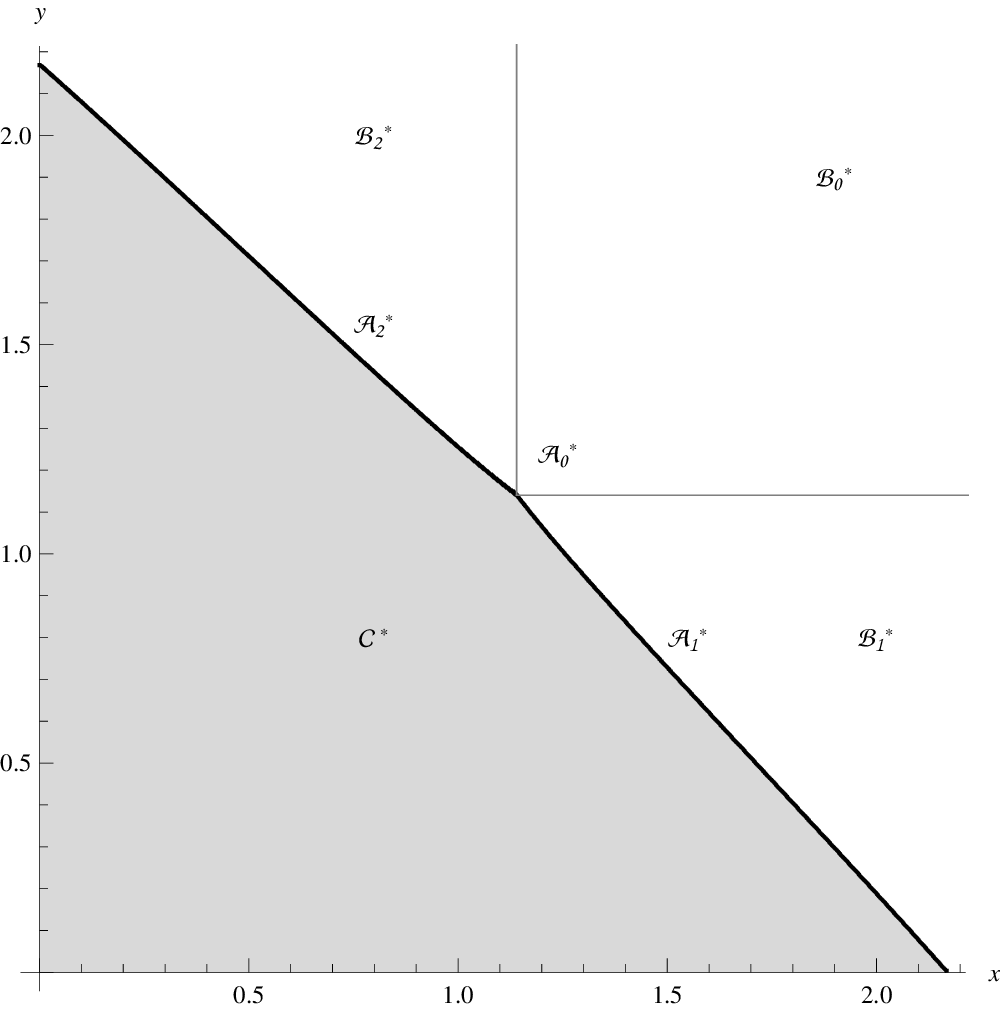}%
\\
Figure 9.1: Optimal Curve Strategy
\end{center}}}%
%EndExpansion
\]

\[%
%TCIMACRO{\FRAME{itbpFU}{5.4172in}{2.84in}{0in}{\Qcb{Figure 9.2: $%
%V(x,y)-\frac{x+y}{2}$ vs. $V_S(x,y)-\frac{x+y}{2}$ vs. $\frac{V_M(x+y)}%
%{2}-\frac{x+y}{2}$}}{}{voptimopapernew.eps}%
%width 5.4172in;  height 2.84in;  depth 0in;  original-width 5.361in;
%original-height 2.7968in;  cropleft "0";  croptop "1";  cropright "1";
%cropbottom "0";  filename 'VOptimoPaperNew.eps';file-properties "XNPEU";}}}%
%BeginExpansion
% {\parbox[b]{5.4172in}{\begin{center}
% \includegraphics[
% height=2.84in,
% width=5.4172in
{\parbox[b]{4in}{\begin{center}
\includegraphics[
height=2.1in,
width=4in
]%
{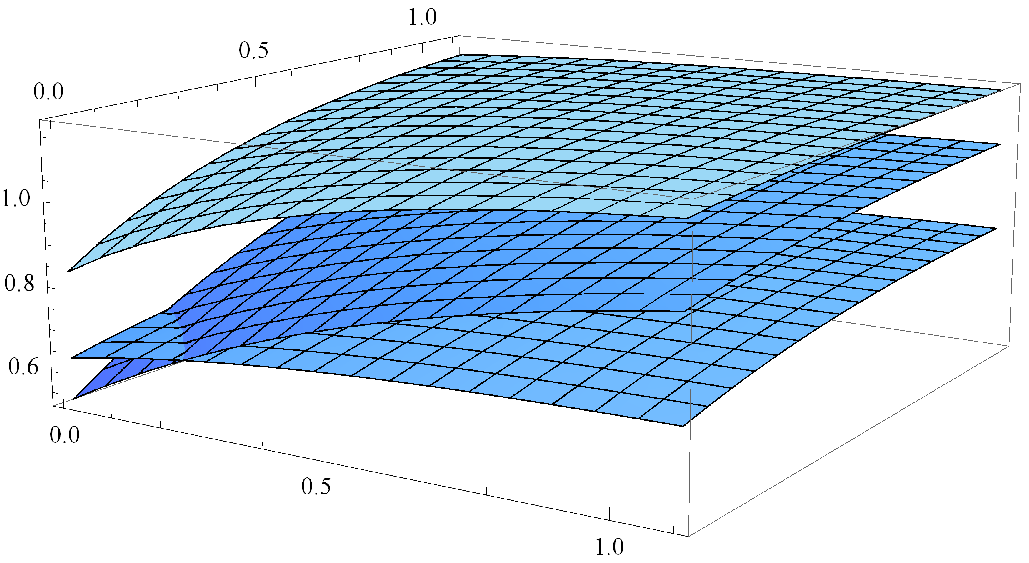}%
\\
Figure 9.2: $V(x,y)-\frac{x+y}{2}$ vs. $V_S(x,y)-\frac{x+y}{2}$ vs. $%
\frac{V_M(x+y)}{2}-\frac{x+y}{2}$
\end{center}}}%
%EndExpansion
\]

\bibliographystyle{abbrv}
\bibliography{ref2c.bib}

\begin{thebibliography}{10}

\bibitem{AlLau}
H.~Albrecher and V.~Lautscham.
\newblock Dividends and the time of ruin under barrier strategies with a
  capital-exchange agreement.
\newblock {\em Preprint, University of Lausanne}, 2014.

\bibitem{AlTho09}
H.~Albrecher and S.~Thonhauser.
\newblock Optimality results for dividend problems in insurance.
\newblock {\em Rev. R. Acad. Cienc. Exactas F\'\i s. Nat. Ser. A Math. RACSAM},
  103(2):295--320, 2009.

\bibitem{asal}
S.~Asmussen and H.~Albrecher.
\newblock {\em Ruin probabilities}.
\newblock Advanced Series on Statistical Science \& Applied Probability, 14.
  World Scientific Publishing Co. Pte. Ltd., Hackensack, NJ, second edition,
  2010.

\bibitem{Avanzi}
B.~Avanzi.
\newblock Strategies for dividend distribution: a review.
\newblock {\em N. Am. Actuar. J.}, 13(2):217--251, 2009.

\bibitem{Avram08b}
F.~Avram, Z.~Palmowski, and M.~Pistorius.
\newblock A two-dimensional ruin problem on the positive quadrant.
\newblock {\em Insurance Math. Econom.}, 42(1):227--234, 2008.

\bibitem{Avram08a}
F.~Avram, Z.~Palmowski, and M.~R. Pistorius.
\newblock Exit problem of a two-dimensional risk process from the quadrant:
  exact and asymptotic results.
\newblock {\em Ann. Appl. Probab.}, 18(6):2421--2449, 2008.

\bibitem{azmu05}
P.~Azcue and N.~Muler.
\newblock Optimal reinsurance and dividend distribution policies in the
  {C}ram\'er-{L}undberg model.
\newblock {\em Math. Finance}, 15(2):261--308, 2005.

\bibitem{azmu2}
P.~Azcue and N.~Muler.
\newblock Minimizing the ruin probability allowing investments in two assets: a
  two-dimensional problem.
\newblock {\em Math. Methods Oper. Res.}, 77(2):177--206, 2013.

\bibitem{azmu}
P.~Azcue and N.~Muler.
\newblock {\em Stochastic Optimization in Insurance: a Dynamic Programming
  Approach}.
\newblock Springer Briefs in Quantitative Finance. Springer, 2014.

\bibitem{BGL}
A.~Badescu, L.~Gong, and S.~Lin.
\newblock Optimal capital allocations for a bivariate risk process under a risk
  sharing strategy.
\newblock {\em Preprint, University of Toronto}, 2015.

\bibitem{badila}
S.~Badila, O.~Boxma, and J.~Resing.
\newblock Two parallel insurance lines with simultaneous arrivals and risks
  correlated with inter-arrival times.
\newblock {\em Insurance Math. Econom.}, 61:48--61, 2015.

\bibitem{cralio}
M.~G. Crandall and P.-L. Lions.
\newblock Viscosity solutions of {H}amilton-{J}acobi equations.
\newblock {\em Trans. Amer. Math. Soc.}, 277(1):1--42, 1983.

\bibitem{czpal}
I.~Czarna and Z.~Palmowski.
\newblock {D}e {F}inetti's dividend problem and impulse control for a
  two-dimensional insurance risk process.
\newblock {\em Preprint,arXiv:0906.2100v3}, 2011.

\bibitem{DeFin57}
B.~De~Finetti.
\newblock Su un' impostazione alternativa dell teoria collettiva del rischio.
\newblock {\em Transactions of the XVth congress of actuaries}, (II):433--443,
  1957.

\bibitem{Ger68}
H.~U. Gerber.
\newblock Entscheidungskriterien fuer den zusammengesetzten
  {P}oisson-{P}rozess.
\newblock {\em Schweiz. Aktuarver. Mitt.}, (1):185--227, 1969.

\bibitem{gershiu06}
H.~U. Gerber and E.~S.~W. Shiu.
\newblock On the merger of two companies.
\newblock {\em N. Am. Actuar. J.}, 10(3):60--67, 2006.

\bibitem{IvBox}
J.~Ivanovs and O.~Boxma.
\newblock A bivariate risk model with mutual deficit coverage.
\newblock {\em Preprint, University of Lausanne}, 2015.
\newblock \texttt{arXiv:1501.02927}.

\bibitem{KuSch}
N.~Kulenko and H.~Schmidli.
\newblock Optimal dividend strategies in a {C}ram\'er-{L}undberg model with
  capital injections.
\newblock {\em Insurance Math. Econom.}, 43(2):270--278, 2008.

\bibitem{loeren}
R.~L. Loeffen and J.-F. Renaud.
\newblock De {F}inetti's optimal dividends problem with an affine penalty
  function at ruin.
\newblock {\em Insurance Math. Econom.}, 46(1):98--108, 2010.

\bibitem{radsh}
R.~Radner and L.~Shepp.
\newblock Risk vs. profit potential: a model for corporate strategy.
\newblock {\em Journal of Economic Dynamics and Control}, 20:1373--1393, 1996.

\bibitem{schmidli08}
H.~Schmidli.
\newblock {\em Stochastic control in insurance}.
\newblock Probability and its Applications (New York). Springer-Verlag London,
  Ltd., London, 2008.

\bibitem{soner}
H.~M. Soner.
\newblock Optimal control with state-space constraint. {I}.
\newblock {\em SIAM J. Control Optim.}, 24(3):552--561, 1986.

\bibitem{ThAl}
S.~Thonhauser and H.~Albrecher.
\newblock Dividend maximization under consideration of the time value of ruin.
\newblock {\em Insurance Math. Econom.}, 41(1):163--184, 2007.

\end{thebibliography}

\end{document}